\documentclass[11pt,reqno,oneside]{amsart}

\usepackage[utf8]{inputenc}
\usepackage[T1]{fontenc}
\usepackage{graphics}
\usepackage[colorlinks=true, pdfstartview=FitV,linkcolor=blue,citecolor=blue,urlcolor=blue]{hyperref}
\usepackage{cleveref}

\usepackage{amsmath,amsthm,amssymb}

\usepackage{fullpage}
\usepackage{numprint}
\npthousandsep{,}

\usepackage{graphicx}
\usepackage{subcaption}
\usepackage[dvipsnames]{xcolor}

\usepackage{booktabs}

\usepackage{tikz}
\usetikzlibrary{arrows.meta, positioning, calc, fit}
\usepackage{graphics}
\usepackage{pgf}

\usepackage{mathtools}
\mathtoolsset{showonlyrefs}

\newenvironment{red}{\relax\color{red}}{\hspace*{.5ex}\relax}
\newenvironment{blue}{\relax\color{blue}}{\hspace*{.5ex}\relax}
\newcommand{\ber}{\begin{red}}
\newcommand{\er}{\end{red}}
\newcommand{\beb}{\begin{blue}}
\newcommand{\eb}{\end{blue}}

\usepackage{enumitem}
\newlist{todolist}{itemize}{2}
\setlist[todolist]{label=$\square$}
\usepackage{pifont}
\theoremstyle{plain}
\newtheorem*{theorem*}{Theorem}

\newtheorem*{proposition*}{Proposition}

\newtheorem*{corollary*}{Corollary}

\theoremstyle{definition}

\DeclareMathOperator{\MCC}{MCC}
\DeclareMathOperator{\ACC}{ACC}

\newcommand{\Q}{\mathbb{Q}}
\newcommand{\F}{\mathbb{F}}

\title{Learning Euler Factors of Elliptic Curves}

\author[Babei]{Angelica Babei}
\address{
  Department of Mathematics,
  Howard University,
  204 Academic Service Building B
  Washington, D.C. 20059
  USA
}
\email{\href{mailto:angelica.babei@howard.edu}{angelica.babei@howard.edu}}
\urladdr{\url{https://angelicababei.com}}

\author[Charton]{Fran\c cois Charton}
\address{
  FAIR, Meta -- CERMICS, Ecole des Ponts
}
\email{\href{mailto:fcharton@meta.com}{fcharton@meta.com}}

\author[Costa]{Edgar Costa}
\address{
  Department of Mathematics,
  Massachusetts Institute of Technology,
  77 Massachusetts Ave.,
  Cambridge,
  MA 02139,
  USA
}
\email{\href{mailto:edgarc@mit.edu}{edgarc@mit.edu}}
\urladdr{\url{https://edgarcosta.org}}

\author[Huang]{Xiaoyu Huang}
\address{
  Department of Mathematics,
  Temple University, Philadelphia, PA, 19122, USA.
}
\email{\href{mailto:xiaoyu.huang@temple.edu}{xiaoyu.huang@temple.edu}}

\author[Lee]{Kyu-Hwan Lee}
\address{Department of Mathematics, University of Connecticut, Storrs, CT 06269, USA \hfill \break \indent Korea Institute for Advanced Study, Seoul 02455, Republic of Korea}
\email{\href{mailto:khlee@math.uconn.edu}{khlee@math.uconn.edu}}
\urladdr{\url{https://khlee-math.github.io/}}

\author[Lowry-Duda]{David Lowry-Duda}
\address{ICERM, 121 S Main St, Providence, RI 02903}
\email{\href{mailto:david@lowryduda.com}{david@lowryduda.com}}
\urladdr{\url{https://davidlowryduda.com}}

\author[Narayanan]{Ashvni Narayanan}
\address{
 Sydney Mathematical Research Institute, University of Sydney, City Road, Camperdown, Sydney, NSW 2006, Australia
}
\email{\href{mailto:ashvni.narayanan@sydney.edu.au}{ashvni.narayanan@sydney.edu.au}}

\author[Pozdnyakov]{Alexey Pozdnyakov}
\address{
  Department of Mathematics,
  Princeton University, Princeton, NJ,  08544-1000, USA.
}
\email{\href{mailto:ap5763@princeton.edu}{ap5763@princeton.edu}}

\date{\today}

\begin{document}

\begin{abstract}
    We apply transformer models and feedforward neural networks to predict
    Frobenius traces $a_p$ from elliptic curves given other traces $a_q$.
    We train further models to predict $a_p \bmod 2$ from $a_q \bmod 2$, and
    cross-analysis such as $a_p \bmod 2$ from $a_q$.
    Our experiments reveal that these models achieve high accuracy, even in the
    absence of explicit number-theoretic tools like functional equations of
    $L$-functions. We also present partial interpretability findings.
\end{abstract}

\maketitle

\section{Introduction}

When must different elliptic curves give different Euler factors (or Frobenius traces $a_p$)?
By a result of Serre, see \cite[Th\'eor\`eme 5]{serre1981quelques} and \cite[Corollary 4.8]{bucur2016application}, under the Riemann hypothesis for Artin $L$-functions, if $E_1$ and $E_2$ are two $\overline{\mathbb{Q}}$-nonisogenous elliptic curves over some number field and neither has complex multiplication, then there exists a prime $p \nmid N$ of size $O(\log^2(N))$ such that $a_p(E_1) \ne a_p(E_2)$, where $N$ is the product of the absolute conductors of $E_1$ and $E_2$.
Thus a finite set of $a_p(E)$ values determines all others.
However, even when enough $a_p(E)$ values are provided to uniquely identify the isogeny class, no efficient algorithm is known for determining the remaining $a_p(E)$ values from this finite set, given a bound on the conductors.

Motivated by this challenge, we investigate the extent to which machine learning (ML) can predict Frobenius traces associated with elliptic curves.
More precisely, given a finite subset $a_p(E)$ values, we explore whether ML models can successfully infer other values and, in doing so, shed light on underlying structural properties of elliptic curves.
The interplay between machine learning and number theory has rapidly expanded in recent years (\cite{AHLOS,HLOa,HLOb,HLOc,HLOP,KV22,Poz,babei2024machine}), and we aim to contribute to this emerging field.
In particular, we will apply transformer models developed by the second-named author~\cite{Char}, along with feed-forward neural networks.

Our experimental approach in \Cref{p: ap} centers on training transformer models to predict $a_p$ for a specific $p$ from a selection of $a_q$ values.
Initial results show that models achieve surprisingly high accuracy when predicting $a_{97}$ from $a_q$ with $q < 97$,
suggesting that some intrinsic patterns in the Frobenius traces are being captured.
We also trained models to predict $a_2$ and $a_3$ for curves with good reduction at these primes.
To better understand what ML models \emph{learn}, we examine principal component analysis (PCA) plots for embeddings in the learning process.
These plots reveal modular behavior, suggesting that the models might be learning congruences modulo $2$ or $4$.
This observation naturally leads us to conduct experiments to learn $a_p \pmod 2$ from $a_q$ values, where we again observe high accuracy.

Building on these observations,
it is natural to try to predict \( a_p \pmod{2} \) from \( \{ a_q \pmod{2} \} \).
We do this in \Cref{p: apmod2}.
However, this approach introduces certain complications.
Notably, reduction modulo 2 produces duplicates in the dataset, causing significant overlap between test and training sets.
Additionally, this reduction results in indeterminacy of \( a_p \pmod{2} \) for identical tuples \( \{ a_q \pmod{2} \} \).
To address these issues, we adopt two strategies: first, we remove duplicates. Second, since the proportion of indeterminate cases is small for conductors up to $10^7$, we retain these cases and conduct experiments with caution.

We also conduct experiments with feedforward neural networks to predict $a_p$
or $(a_p \bmod 2)$
from either $\{ a_q : q \neq p, q < 100 \}$, $\{ a_q \bmod 2 : q \neq p, q < 100 \}$ or $\{ \frac{a_q}{\sqrt{q}} \bmod 2 : q \neq p, q < 100 \}$.
Models have comparable or lower performance.
Saliency analysis suggests that the most important features might depend
on input encoding.
We describe these experiments and analysis further in \Cref{p: NN}.

It is surprising that ML models often successfully predict $a_p \pmod \ell$ from
remaining $a_q \pmod \ell$ values, as in this case we don't have typical tools like $L$-functions and functional equations to directly solve for missing coefficients.
Instead, one needs to enumerate all possible mod $\ell$ Galois representations to deduce the right answer.

Our work provides partial but compelling experimental evidence that ML can
effectively learn relationships among Euler factors and Frobenius traces.  We provide evidence that our transformer models' predictions of \( a_p \) rely on implicitly predicting \( a_p \pmod{2} \) values.
However, given the opaque nature of neural networks and transformers, the
precise mechanisms underlying these predictions remain an open mystery at the
moment and calls for further analysis.

\section*{Acknowledgments}

We thank the Harvard Center of Mathematical Sciences and Applications and to the organizers of the Mathematics and Machine Learning Program in fall of 2024 for creating such an exciting, welcoming environment.
We would also like to thank Mike Douglas and Drew Sutherland for their support, guidance, and advice.

Costa was supported by Simons Foundation grant SFI-MPS-Infrastructure-00008651 and by Simons Foundation grant 550033.
KHL was partially supported by a grant from the Simons Foundation (\#712100).
DLD was supported by the Simons Collaboration in Arithmetic Geometry, Number
Theory, and Computation via the Simons Foundation grant 546235. Ashvni thanks the Hooper-Shaw foundation and the Sydney Mathematical Research Institute for supporting her.

This research includes calculations carried out on HPC resources supported in part by the National Science Foundation through major research instrumentation grant number 1625061 and by the US Army Research Laboratory under contract number W911NF-16-2-0189.

\section{ML architectures and preliminaries}
\label{p: prelims}

\subsection{Elliptic curves}
An elliptic curve over $\Q$ is a smooth, projective curve $E$ of genus 1 with a distinguished point at infinity $O \in E(\Q)$. It  can be described by a Weierstrass equation:
\[ E: y^2 = x^3+Ax+B, \qquad A, B \in \Q.\]

One can associate to $E$ an $L$-function, given by the Euler product
\[L(E, s) = \prod_p L_p(E, s). \]

At all but finitely many primes $p$, where $E$ has good reduction,  the local factor is
\[ L_p(E, s) = \frac{1}{1-a_p(E)p^{-s}+p^{1-2s}} \]
with the trace of Frobenius at $p$ equalt to \[ a_p(E) = p+1-| E(\F_p)|,\]
and  $|E(\F_p)|$ is the number of solutions to the defining equation  modulo $p$ . At the remaining finitely many primes, $E$ has bad reduction.
The conductor $N(E)$ is a product of powers of these primes with bad reduction, with exponents being dictated by whether $E$ has multiplicative or additive reduction.

In the remaining sections, we perform experiments to predict the values $a_p$, which determine the Euler factors, given traces of Frobenius $a_q$ for nearby primes $q \ne p$.

\subsection{Model Architecture: Transformers}
\label{s: models}

The transformer architecture, introduced in \emph{Attention is All You Need} \cite{vaswani2017attention}, has revolutionized deep learning by eliminating recurrence and instead leveraging self-attention mechanisms to process sequences in parallel. This design has proven highly effective in natural language processing (NLP) and is increasingly being applied across diverse domains. While we refer the reader to, for example,  \cite{elhage2021mathematical} for details on transformers, we briefly overview their architecture below.

A transformer consists of several key components:
\begin{enumerate}
\item {\bf token embedding layer}:
The first layer maps discrete tokens to continuous vector representations in a high-dimensional space. Given a vocabulary of size $v$ and an embedding dimension $d$, each token is represented as a $d$-dimensional vector by multiplying its one-hot representation by a learnable $v \times d$ matrix.

\item {\bf positional encoding}:
Since transformers lack recurrence, positional encodings are added to embeddings to introduce sequential information.
These embeddings are learned, just like the token embedding.

\item {\bf self-attention mechanism}:
Self-attention enables each token to attend to all others in the sequence, allowing the model to capture long-range dependencies. The attention scores are computed using:
\begin{equation}
    \text{Attention}(Q, K, V) = \text{softmax}\left( \frac{QK^T}{\sqrt{d_k}} \right)V
\end{equation}
where $Q$, $K$, and $V$ represent queries, keys, and values derived from input embeddings and $d_k$ is the dimension of the key vectors. Here, the softmax function ensures that attention scores form a probability distribution.

\item {\bf multiple attention heads}:
Transformers employ multiple attention heads, each focusing on different aspects of the input. These heads operate independently, and their contributions are summed before passing through subsequent layers. Given $h$ attention heads, the multi-head attention mechanism is expressed as:
\begin{equation}
    \text{MultiHead}(Q, K, V) = \sum_{i=1}^{h} \text{Attention}(Q W_i^Q, K W_i^K, V W_i^V) W_i^O
\end{equation}
where $W_i^Q, W_i^K, W_i^V,$ and $W_i^O$ are learned projection matrices for the $i$-th head.  Attention heads can be interpreted as information-moving operations, transferring data between tokens. This movement is crucial for understanding how transformers capture relationships across sequences.

\item {\bf feedforward network}:
Each attention layer is followed by a position-wise feedforward network (FFN) containing one hidden layer of dimension $4d$,
applied independently to each sequence position.

\item {\bf residual connections and layer normalization}:
Residual connections ensure stable gradient flow, while layer normalization improves convergence.

\end{enumerate}

The transformer processes input sequences by first embedding the tokens and
positions as two $d$-dimensional vectors, which are added together and fed into
$n$ transformer layers, composed of a multi-head attention and a feed-forward
network.
The output of the last layer is then decoded by a linear layer similar to the
token encoder, which predicts, for each position in the output, the
probabilities of the corresponding output token.
The most likely output is the model prediction.

\subsection{Int2Int}
\label{ss: int2int}
We utilize Int2Int \cite{Char}, which implements a sequence-to-sequence transformer that is specialized to ``translate'' short sequences of integers into short sequences of integers. Our model is configured as an encoder-only architecture, as shown in \Cref{fig:archi}, with a total of \numprint{4223805} trainable parameters. We adopt the same encoding as the \textit{PositionalInts} encoding in Int2Int for the inputs, where each integer is represented by a token indicating its sign followed by a token for its absolute value. For instance, the list of integers $-1, 2, 5$ is converted to the string ``- 1 + 2 + 5''. To handle output values, we apply a right shift so that all values are nonnegative. For instance, the smallest output value $-2$ is mapped to the string ``0'' in the $a_2$ prediction task.

\definecolor{DarkGray}{HTML}{68904d}
\definecolor{Orange}{HTML}{ffae42}
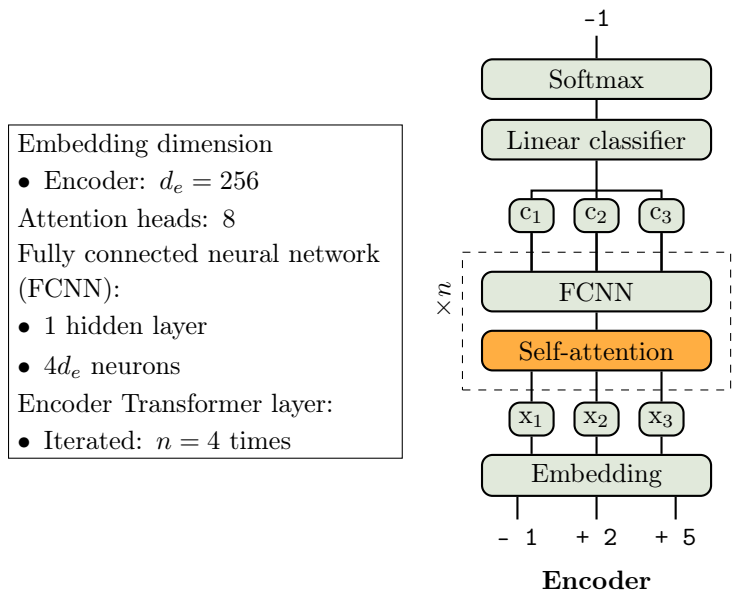
\begin{figure}[ht]
\small
\centering
\begin{tikzpicture}[
    box/.style={rectangle,draw,fill=DarkGray!20,node distance=1cm,text width=8em,text centered,rounded corners,minimum height=1.5em,thick},
    box2/.style={rectangle,draw,fill=DarkGray!20,node distance=1cm,text width=1em,text centered,rounded corners,minimum height=0.5em,thick},
    % arrow/.style={draw,-latex',thick},
  ]
  \tikzset{arrow/.style={-latex, thick}}
  \node [box] (ffne) {FCNN};

  \node [box2,above=0.5 of ffne] (tok2) {c$_2$};
  \node [box2,left=0.25 of tok2] (tok1) {c$_1$};
  \node [box2,right=0.25 of tok2] (tok3) {c$_3$};
  % \node [above=0.2 of tok2]{Encoder output};
  \node [box,above=0.5 of tok2] (linear){Linear classifier};
  \node [box,above=0.25 of linear] (softmax){Softmax};
  \node [above=0.3 of softmax] (output) {\texttt{-1}};

  \path[draw,thick] (linear) -- (softmax);
  \path[draw,thick] (softmax) -- (output);
  \path[draw,thick] (tok1.south) -- (ffne.north -| tok1.north);
  \path[draw,thick] (tok2.south) -- (ffne);
  \path[draw,thick] (tok3.south) -- (ffne.north -| tok3.north);

  \node [box,below=0.23 of ffne,fill=Orange] (selfatt){Self-attention};  
  \node [label={[rotate=90]160:$\times n$ },rectangle,draw,dashed,inner sep=0.7em,fit=(selfatt) (ffne)] (encoder1) {};

  \node [box2,below=0.4 of selfatt] (se){x$_2$};
  \node [box2,left=0.25 of se] (se2){x$_1$};
  \node [box2,right=0.25 of se] (se3){x$_3$};
  \node [box,below=0.23 of se] (embed) {Embedding};
  \node [below=0.3 of embed] (inp2) {\texttt{+ 2}};
  \node [right=0.25 of inp2] (inp3) {\texttt{+ 5}};
  \node [left=0.25 of inp2] (inp1) {\texttt{- 1}};
  \path[draw,thick] (inp2) -- (embed);
  \path[draw,thick] (inp3) -- (embed.south -| inp3);
  \path[draw,thick] (inp1) -- (embed.south -| inp1);
  \node [below=0.1 of inp2]{\textbf{Encoder}};

  % \coordinate[above right =0.2 and 0.7 of tok3](conn);
  \coordinate[above =0.1 of tok2](conn);
  \path[draw,thick] (tok1) |- (conn);
  \path[draw,thick] (tok2) |- (linear.south);
  \path[draw,thick] (tok3) |- (conn);

  \path [draw,thick] (ffne.north -| tok1) -- (tok1);
  \path [draw,thick] (ffne) -- (tok2);
  \path [draw,thick] (ffne.north -| tok3) -- (tok3);
  \path [draw,thick] (selfatt) -- (ffne);
  \path [draw,thick] (se) -- (selfatt);
  \path [draw,thick] (se2) -- (selfatt.south -| se2);
  \path [draw,thick] (se3) -- (selfatt.south -| se3);
  \path[draw,thick] (embed.north -| se2) -- (se2);
  \path[draw,thick] (embed) -- (se);
  \path[draw,thick] (embed.north -| se3) -- (se3);

  \node[draw,text width=5cm] at (-5.2,0) {Embedding dimension\\
  \begin{itemize}[leftmargin=*, itemsep=1pt, topsep=1pt, parsep=1pt]
  \item Encoder: $d_e=256$ \\
  \end{itemize}
  Attention heads: $8$

  Fully connected neural network (FCNN):\\ 
  \begin{itemize}[leftmargin=*, itemsep=1pt, topsep=1pt, parsep=1pt]
  \item 1 hidden layer \\
  \item $4d_e$ neurons \\
  \end{itemize}
  Encoder Transformer layer:\\
  \begin{itemize}[leftmargin=*, itemsep=1pt, topsep=1pt, parsep=1pt]
      \item Iterated: $n = 4$ times
  \end{itemize}
  };
\end{tikzpicture}
\vspace{-0.1cm}
\caption{We implement an encoder only transformer architecture using Int2Int.}\label{fig:archi}
%\vspace{-0.1cm}
\end{figure}

\subsection{Performance metrics}
\label{s: performance}

\subsubsection{Accuracy}
\label{s: acc}

In the case of binary classification on balanced datasets, the most instructive performance metric is accuracy, which measures the proportion of correct predictions. We use this metric in the experiments in Sections \ref{p: apmod2} and \ref{p: NN}.

\subsubsection{Matthews Correlation Coefficient}
\label{s: mcc}

Another performance metric is the the Mathews correlation coefficient (MCC), which is especially useful for unbalanced binary classification experiments. We use the MCC and a modified version as described below in all our experiments.

For binary classification experiments, the MCC is defined as in \cite{matthews1975comparison} by
$$\MCC := \frac{TP \times TN - FP \times FN}{\sqrt{(TP + FP)(TP + FN)(TN + FP)(TN + FN)}},$$
where $TP, FP, TN$ and $FN$  are the number of true positives, false positives, true negatives and false negatives, respectively. The MCC attains a value in the range $[-1, 1]$, where $1$ indicates perfect prediction,  $0$ indicates no correlation between predictions and true values, and $-1$ indicates perfectly wrong classification. For a balanced binary classification experiment, one obtains the accuracy ACC is obtained by $\ACC = \frac{\MCC + 1}{2}$.

One of the generalizations of MCC to multiclass experiments is given in \cite{gorodkin2004comparing} by the formula

$$\MCC := \frac{\sum_k \sum_l \sum_m C_{kk}C_{lm} - C_{kl}C_{mk}}{\sqrt{\sum_k
 (\sum_l C_{kl}) (\sum_{k'| k' \ne k} \sum_{l'} C_{k'l'})} \sqrt{\sum_k
 (\sum_l C_{lk}) (\sum_{k'| k' \ne k} \sum_{l'} C_{l'k'})} },$$
 where $C_{ij}$ is the $(i,j)$ entry in the confusion matrix $C$. This generalization is implemented in sklearn, and is the one we use in our multiclass experiments.

\subsection{Datasets}
\label{s: datasets}

We concentrated our efforts on isogeny classes of elliptic curves with small conductors, specifically those less than \(10^7\). In particular, we started with a dataset created by Andrew Sutherland that extends the one in \cite{stein-watkins}, constructed using similar methods and including elliptic curves with conductors up to \(10^8\). From this superset, we extracted two subsets, \texttt{ECQ6} and \texttt{ECQ7} \cite{ECQ7}, of sizes \numprint{3607212} and \numprint{28257805}, respectively, corresponding to the isogeny classes with conductors less than \(10^6\) and \(10^7\).

Some experiments in \Cref{p: NN} are performed on isogeny classes of curves with conductor between $2^{15}$ and $2^{16}$ in the LMFDB \cite{lmfdb}. These curves form a subset of \texttt{ECQ6}.

\section{Learning \texorpdfstring{$a_p$}{ap} for small primes}
\label{p: ap}

In this section, we examine if one can learn the true value $a_p$ given nearby traces of Frobenius. We restrict ourselves to an input sequence $(a_q)_{q \ne p, q \le 100}$ of 24 primes and curves of conductor less than $10^6$.

The sequential nature of the data $(a_q)_{q \ne p, q\le 100}$ is well-suited for sequence-to-sequence models such as recurrent neural networks (RNNs) or transformers. With some limitations, transformers have been taught number theory tasks such as modular arithmetic \cite{saxena2024teaching} and predicting the greatest common divisor \cite{charton2023can}. As seen in \cite{charton2023can}, the ability of transformers to learn problem-specific embeddings could be particularly useful for number theory applications.

\subsection{Methodology}

For the following three experiments, we employed a transformer model as implemented in Int2Int ~\cite{Char} with an encoder-only architecture. The vocabulary consists of integers $ 0 \le n \le 19$, signs $+, -$, as well as special tokens. The model architecture follows the description in \Cref{ss: int2int}. We used the Adam optimizer and a learning rate of $3 \cdot 10^{-5}$.  Each task is treated as a classification task, with $39$ classes for predicting $a_{97}$, $5$ classes for $a_2$, and $7$ classes for $a_3$.

The training examples consist of sequences $(a_q)_{q \ne p, q <100}$ of length $24$, for which we predict $a_p$. This data represents isogeny classes of curves with good reduction at $p$ from the dataset \texttt{ECQ6}. In each case, we reserved a test set of 10k isogeny classes sampled uniformly at random, distinct from the training set.

\subsection{Results}

\begin{figure}[ht]
    \centering
        \includegraphics[height=7cm, width=0.5\textwidth]{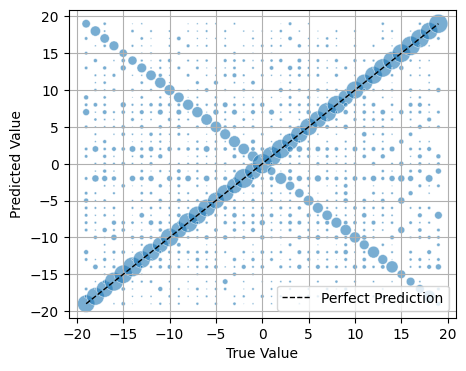}
        \caption{Predicting $a_{97}$. Proportion of cases where the model predicts a given value as a function of its true value.}
    \label{fig: a97_conf}
\end{figure}

We describe in detail the results of the experiment for learning $-19 \le a_{97} \le 19$.
After 100 epochs of 3.5 million training examples, the transformer model reached a maximum accuracy of $0.4923$ and maximum MCC of $0.4711$. In Figure \ref{fig: a97_conf}, we show the proportion of cases where the model predicts a given value as a function of its true value. The main diagonal corresponds to perfect predictions.

Although the experiment for predicting $a_{97}$ performed well compared to the baseline precision of $0.082$, which is the proportion of the most common class, it cannot serve as a prototype for generating sequential values $a_{97}, a_{101}, a_{103}, \dots$. One of the potential causes of low performance could be the large number of classes to be predicted. Therefore, we trained two additional models to predict $a_p$ which contains fewer classes, in particular $a_2$ and $a_3$.

\begin{figure}[ht]

    \centering
    \begin{subfigure}{0.4\textwidth}
        \centering
        \includegraphics[width=\linewidth]{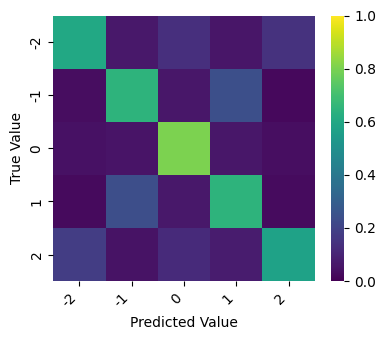}
    \end{subfigure}
    \hspace{0.1\textwidth}
    \begin{subfigure}{0.4\textwidth}
        \centering
        \includegraphics[width=\linewidth]{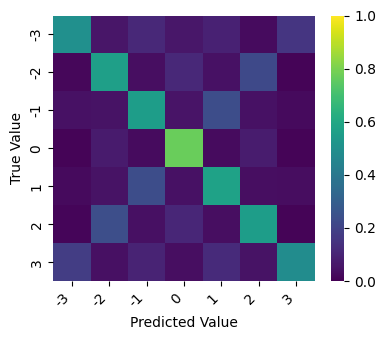}
    \end{subfigure}
    \caption{Confusion matrix heat map for predicting $a_2$ (left) and $a_3$ (right).}
    \label{fig: a2a3_conf}

\end{figure}

The results for predicting $-2 \le a_2 \le 2$, $-3 \le a_3 \le 3$ and $-19 \le a_{97} \le 19$ are summarized in Table \ref{tab: learn_ap}. In addition, Figure \ref{fig: a2a3_conf} shows the heat maps of the normalized confusion matrices for predicting $a_2$ and $a_3$. The alternating color intensity along the rows indicates that the predicted values frequently match the true values modulo 2. We also explore this modular phenomenon further in Section \ref{p: apmod2}, where we predict $a_p \pmod{2}$. Moreover, we visualize the 2D PCA projections of their learned embeddings in \Cref{p: interpret} to understand the underlying modular structure of these embeddings.

We note that in both Figures \ref{fig: a97_conf} and \ref{fig: a2a3_conf}, the models often correctly predict $a_p$ up to sign. In Table \ref{tab: learn_ap}, we also record the sign-agnostic MCC of $|\text{predicted } a_p| \overset{?}{=} |\text{true value of } a_p|$.

In all three cases, the models struggle to distinguish the sign of $a_p$. In \Cref{ss: a2NN}, we perform a binary classification neural network-based experiment to predict the value $a_2 = \pm 1$.

\begin{table}[!ht]
\centering
\begin{tabular}{|c|c|c|c|c|c|c|}
\hline
$p$ & $N_{\text{samples}}$ & $N_{\text{epochs}}$ & Epoch size & $N_{\text{classes}}$ & Max Test MCC & Max Sign-Agnostic MCC\\ \hline
2   & 722,204    & 200 & 700,000    & 5      & 0.5822  & 0.7345 \\ \hline
3   & 1,210,417   & 200 & 1,200,000     & 7      & 0.5205 & 0.6908  \\ \hline
97   & 3,550,727    & 100 & 3,500,000    & 39     & 0.4711 & 0.6266 \\ \hline
\end{tabular}
\caption{Predicting $a_p$ for $p \in \{2, 3, 97\}$ with an encoder-only transformer. We record both the maximum MCC as well as the maximum sign-agnostic MCC.}
\label{tab: learn_ap}
\end{table}

\section{Learning \texorpdfstring{$a_p \pmod{2}$}{ap mod 2}}
\label{p: apmod2}

The confusion matrix results in \Cref{fig: a97_conf} and \Cref{fig: a2a3_conf}  indicate that the models' inference of \( a_p \) rely, at least in part, on implicitly inferring \( a_p \pmod{2} \). In particular, the models should be capable of explicitly determining \( a_p \pmod{2} \). To validate this hypothesis, we conducted two sets of experiments focused on learning \( a_p \pmod{2} \) for all $p < 100$.

\subsection{Given \texorpdfstring{$\{a_q\}_{q \ne p, q < 100}$}{\{aq\}, where p != q < 100,} for curves with good reduction at \texorpdfstring{$p$}{p}}

\begin{figure}[ht]
    \centering
    \includegraphics[height=5.5cm, width=0.48\linewidth]{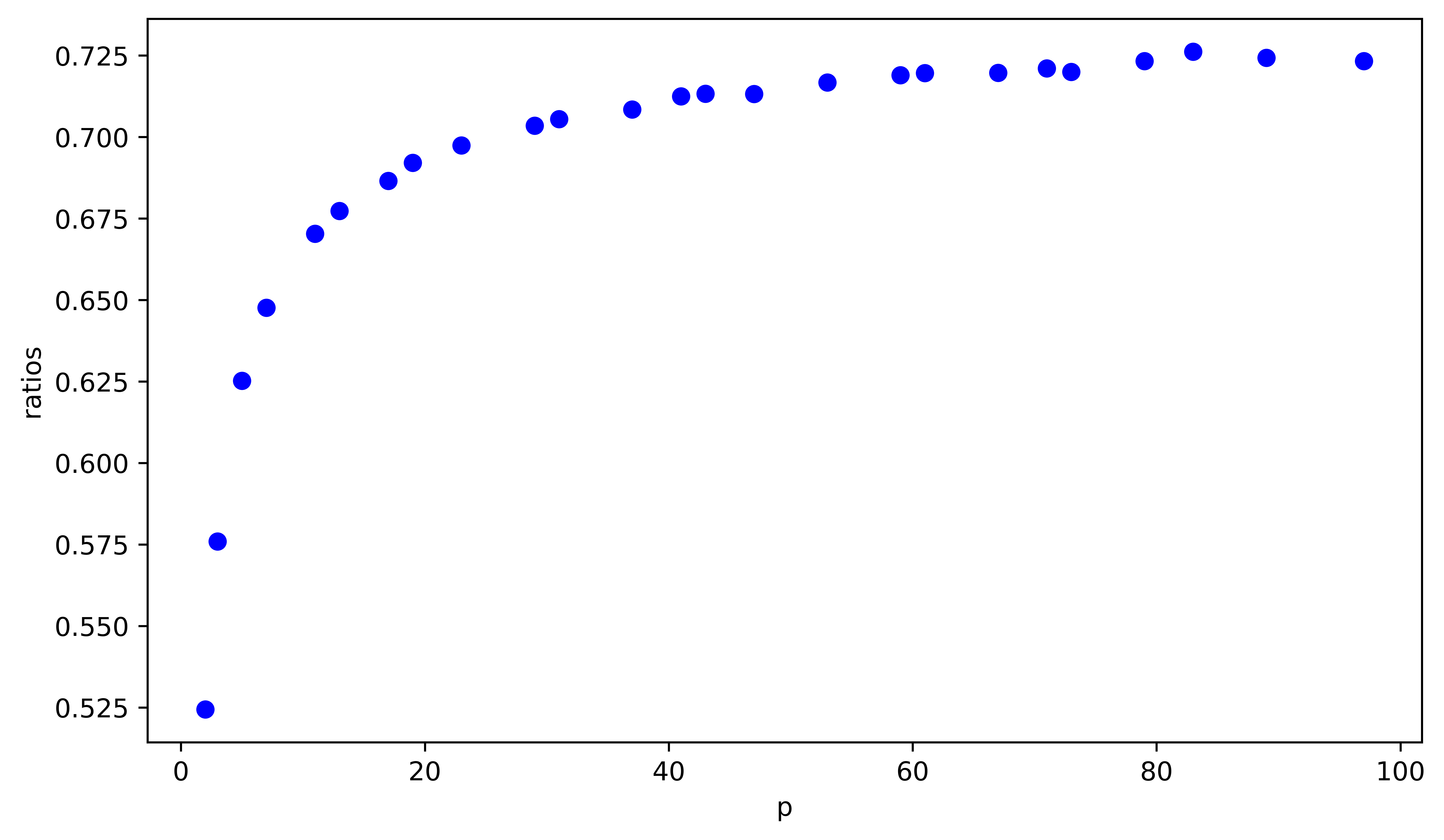}
    \caption{The ratios of curves with $a_p \equiv 0 \pmod 2$ for various $p < 100$ when curves have good reduction at $p$. There is a consistent bias of $a_p \equiv 0 \mod 2$ when $p$ is not ``small''. }
    \label{fig: ap_mod2_ratios}
\end{figure}

In each case of $p$, the training set includes 2 million curves with good reduction at $p$, except for $p = 2, 3, 5$, sourced from \texttt{ECQ6}. The test set consists of 10,000 distinct curves not present in the training set. For $p = 2, 3, 5$, the training set sizes are \numprint{632205}, \numprint{1120417}, and \numprint{1906257}, respectively, due to the good reduction condition. The dataset is intentionally unbalanced for this experiment for two reasons. First, to reflect real-world deployment scenarios where class distributions are naturally imbalanced. Second, the majority class ratio (where $a_p \equiv 0 \mod 2$) ranges from 52.45\% for $p = 2$ to 72.62\% for $p = 97$ in the combined training and test datasets, as shown in \Cref{fig: ap_mod2_ratios}, which is within acceptable limits.

\begin{figure}[ht]
    \centering
    \begin{subfigure}[b]{0.48\linewidth}
        \centering
        \includegraphics[height=5.5cm, width=\linewidth]{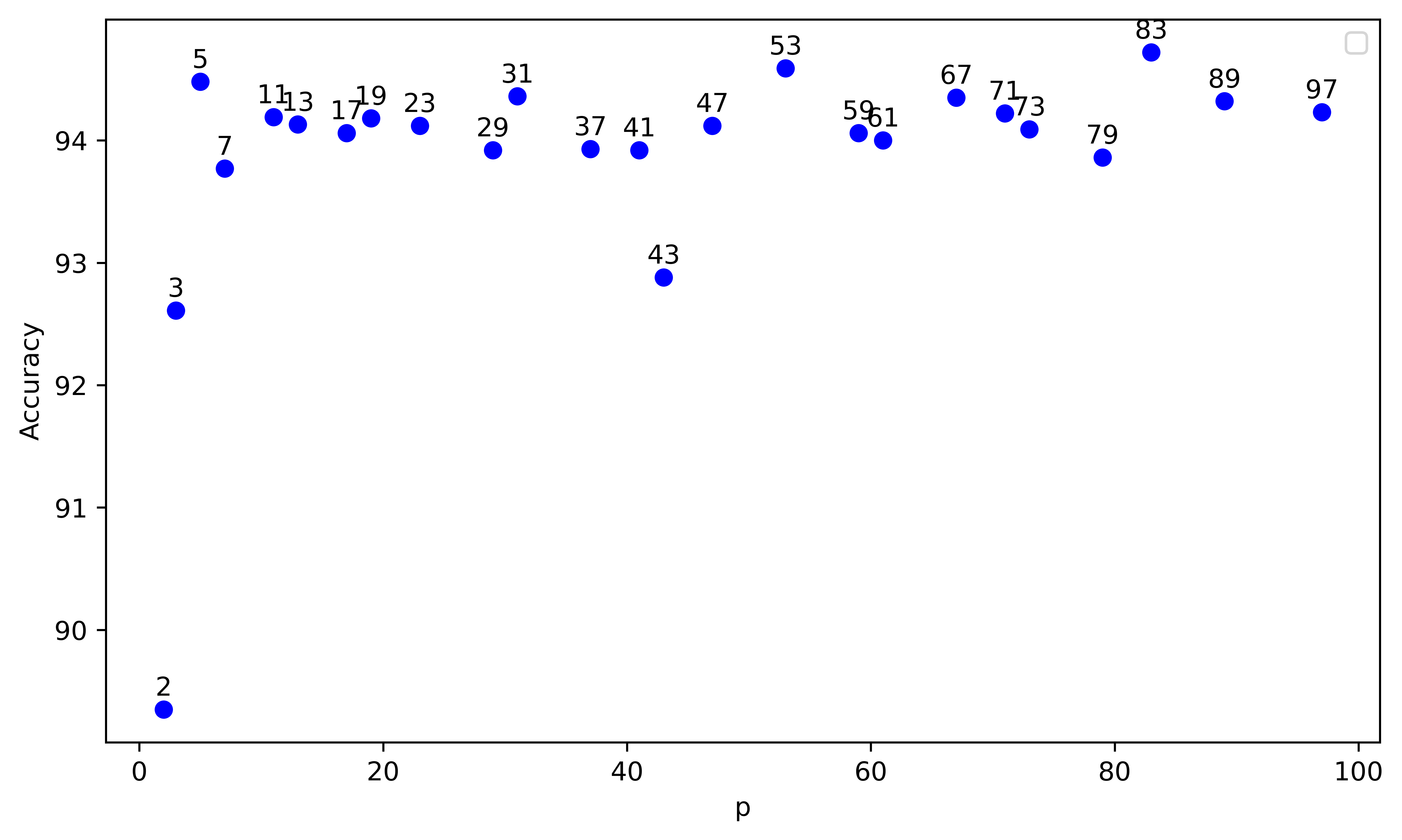}
        \caption{Accuracies}%
    \end{subfigure}
    \hfill
    \begin{subfigure}[b]{0.48\linewidth}
        \centering
        \includegraphics[height=5.5cm,width=\linewidth]{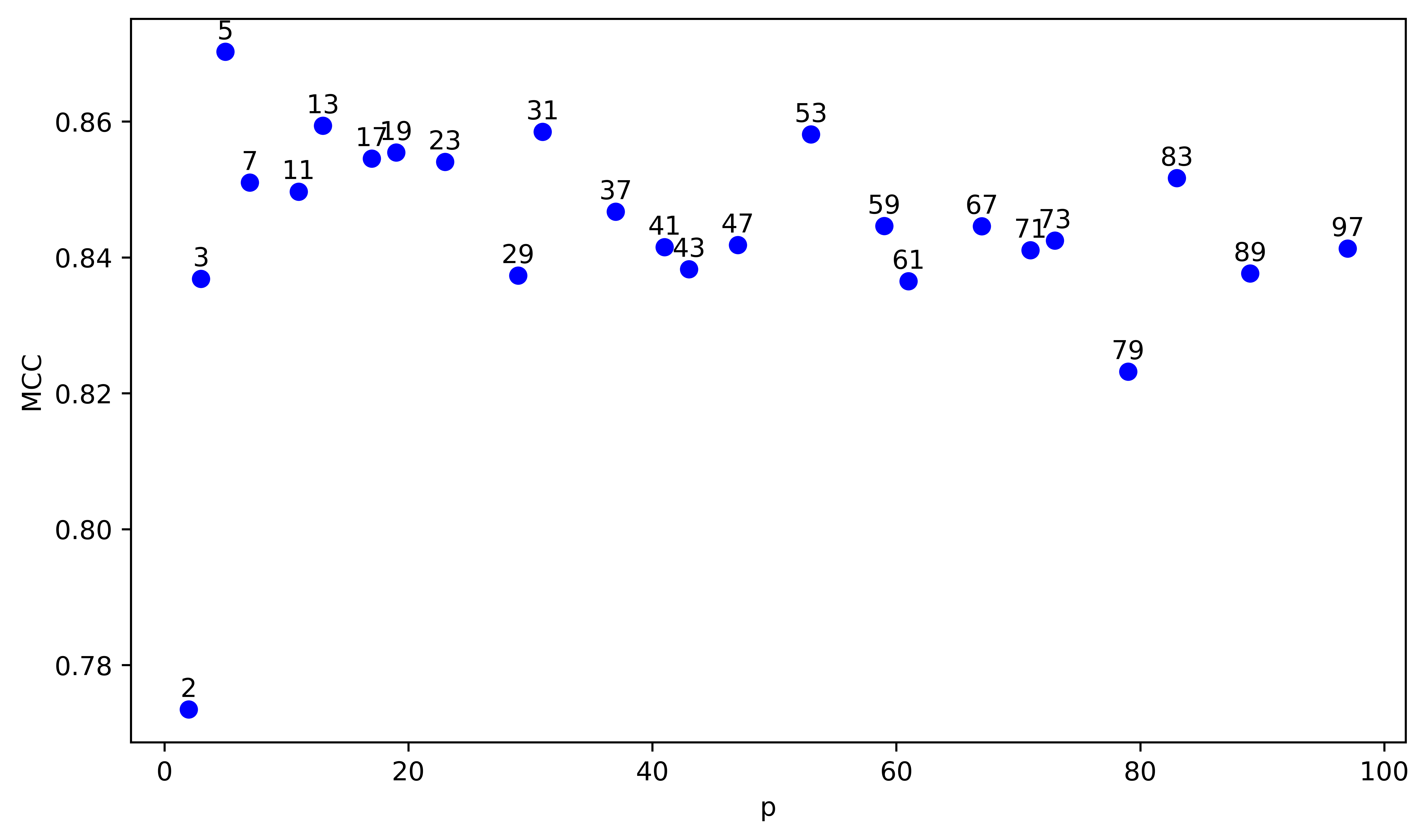}
        \caption{MCCs}
    \end{subfigure}%
    \caption{The results of predicting \( a_p \mod 2 \) using \(\{a_q\}_{q \ne p, q < 100}\) for curves with good reduction at \( p \) show that, except for \( p = 2, 3, 43 \), the model achieves accuracies close to 0.94. The highest accuracy, 0.9472, occurs at \( p = 83 \), while the lowest, 0.8935, is observed at \( p = 2 \). Similarly, the model achieves MCC around 0.84, except for \( p = 2, 79 \). The highest MCC, 5, 0.8703, is recorded at \( p = 5 \), whereas the lowest, 0.7735, occurs at \( p = 2 \).
    }
    \label{fig: raw_ap_res}
\end{figure}

The transformer model was implemented in Int2Int as described in \Cref{ss: int2int}. We used the Adam optimizer and a learning rate of $5 \cdot 10^{-5}$ for the results here. The accuracy and MCC results are presented in \Cref{fig: raw_ap_res}.

Interestingly, the 2D PCA projections of the learned embeddings closely resemble those from \Cref{p: ap}, both presented in \Cref{p: interpret}, where the models predict the true \( a_p \). As discussed in \Cref{p: interpret}, this similarity suggests that the models in \Cref{p: ap} may be inferring \( a_p \) values based on their mod 2 representations.

\subsection{Given \texorpdfstring{$(a_q \pmod{2})_{q\ne p, q<100}$}{\{aq mod 2\}, where p != q < 100,} for curves of good reduction at \texorpdfstring{$p$}{p} and removing duplicates.}
Comparing the previous experiment with one on the same dataset, but whose input is the sequence modulo 2, $(a_p \pmod{2})_{p \ne q, p<100}$), poses two problems.

The first problem is of significant overlap between training and test data. Despite the fact that theoretically there are $\sim 33.5$ million binary tuples of length 25, \texttt{ECQ6} only contains $\sim 0.7$ million representatives modulo $2$. Once we reduce the dataset modulo 2, more than $80\%$ of the test data is seen in the training data. The second issue is one of indeterminacy, as the same input sequence $(a_q\pmod{2} )_{q\ne p, q<100}$  can yield different output $a_p \pmod{2}$.

We create a new dataset that addresses the first issue by removing duplicates, i.e. we select a single representative for a given tuple $(a_q\pmod{2} )_{q<100}$ containing both input and output. However, this significantly reduces the dataset size when restricted to isogeny classes in \texttt{ECQ6}. Therefore, for this series of experiments only, we use curves from \texttt{ECQ7} instead. In particular, we include 5 million isogeny classes of curves, containing all representatives mod 2 in \texttt{ECQ6} as well as approximately $4.3$ million representatives in \texttt{ECQ7} with the lowest conductor.

In the following experiments, we train models to predict $a_p \pmod{2}$ for each prime $p<100$ on the dataset described above. For each prime, we first remove curves with bad reduction at $p$, then balance the dataset according to $a_p \pmod{2}$. After balancing, the training set at each prime has between 5 -- 11.5 \% indeterminates; when restricting to curves of conductor up to $10^6$, the proportion of indeterminates is reduced to 1 -- 2.5 \%.

We reserve \numprint{15000} examples for a test set, distinct from the training set, where \numprint{5000} curves in each of the following conductor ranges:  $10^4 \le N(E) < 10^5$,  $10^5 \le N(E) < 10^6$ and  $10^6 \le N(E) < 10^7$. We use an encoder-only transformer as described in \Cref{ss: int2int}. We ran each experiment for \numprint{1000} epochs of size \numprint{100000} training examples.

\begin{figure}[!ht]
    \centering
        \includegraphics[width=1\textwidth]{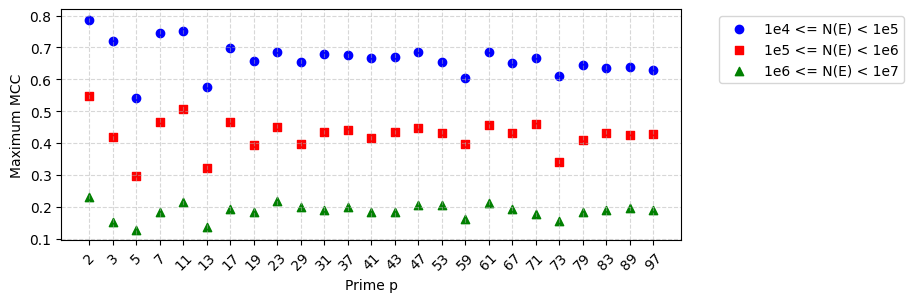}

        \caption{Maximum MCC achieved in three conductor ranges for $a_p \pmod{2}$ for various primes $p < 100$.}
         \label{fig: apmod2accdistinct}
\end{figure}

Figure \ref{fig: apmod2accdistinct} shows the maximum attained MCC across the three conductor ranges for each prime. As expected, the performance is reduced for curves of the highest conductor range $10^6 \le N(E) < 10^7$.

\section{Neural network experiments}
\label{p: NN}

To explore the significance of the transformer architecture used above, we provide several experiments that use a standard feedforward neural network (FFN) to predict $a_p$. For one set of experiments, we use an FFN to predict $a_p \bmod \ell_p$ for $p \in \{2, 97\}$ and $\ell_p \in \{0,2\}$ \footnote{We use the convention $a \equiv b \pmod{0} \iff a=b$.}. As features we use $\{a_q \bmod \ell_q : q < 10^2, q \neq p\}$ where $\ell_q \in \{0,2\}$. The models were trained and tested on all of the curves in \texttt{ECQ6} with good reduction at $p$ and with duplicates removed. All neural networks have layers of width $[2^{n+m},2^{n+m-1}, \cdots, 2^{n+1}, 2^{n}]$ for some $m,n$. All of the features and labels were represented by a one-hot encoding. We always use $N_{\text{samples}} - 10^4$ curves for training and $10^4$ curves for testing. We always use AdamW optimizer with learning rate $10^{-4}$ and weight decay $0.1$, and we train for $10^4$ epochs.

\begin{table}[!ht]
\centering
\begin{tabular}{|c|c|c|c|c|c|c|}
\hline
$p$ & $\ell_p$ & $\ell_q$ & $N_{\text{params}}$ & $N_{\text{samples}}$ & Test Acc & Test MCC \\ \hline
2   & 0        & 0        & \numprint{136581}            & \numprint{732205}             & 0.4136   & 0.2412   \\ \hline
2   & 2        & 0        & \numprint{136482}            & \numprint{732205}             & 0.7364   & 0.4770  \\ \hline
2   & 2        & 2        & \numprint{16674}            & \numprint{292107}             & 0.6651   & 0.3300  \\ \hline
97  & 0        & 0        & \numprint{1709353}            & \numprint{2000000}             & 0.2602   & 0.2268  \\ \hline
97  & 2        & 0        & \numprint{136482}            & \numprint{2000000}            & 0.8691   & 0.6524  \\ \hline
97  & 2        & 2        & \numprint{55842}            & \numprint{697051}             & 0.7433   & 0.4228  \\ \hline
\end{tabular}
\caption{Neural network experiments}
\label{tab:NN_performance}
\end{table}

\begin{figure}[ht]
    \centering
    \begin{subfigure}[b]{0.3\textwidth}
        \centering
        \includegraphics[width=\textwidth]{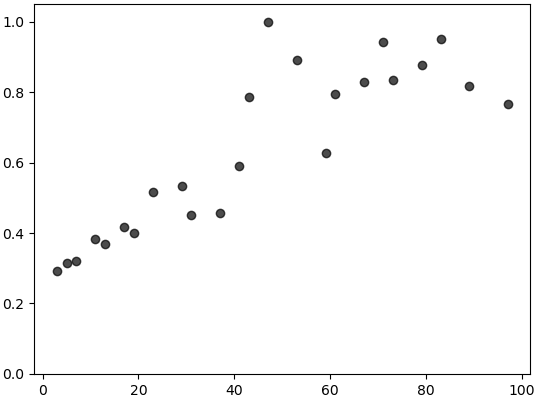}
        \caption{$p=2, \ell_p = 0, \ell_q=0$}
    \end{subfigure}
    \begin{subfigure}[b]{0.3\textwidth}
        \centering
        \includegraphics[width=\textwidth]{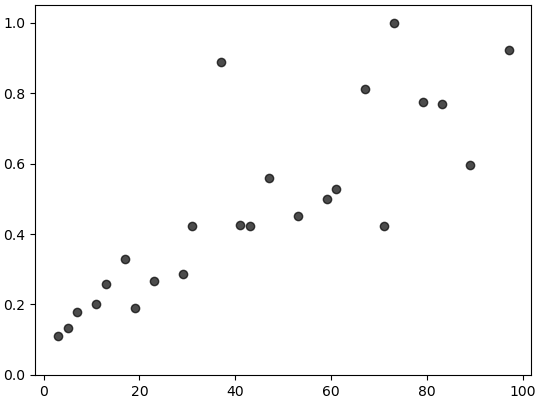}
        \caption{$p=2, \ell_p = 2, \ell_q=0$}
    \end{subfigure}
    \begin{subfigure}[b]{0.3\textwidth}
        \centering
        \includegraphics[width=\textwidth]{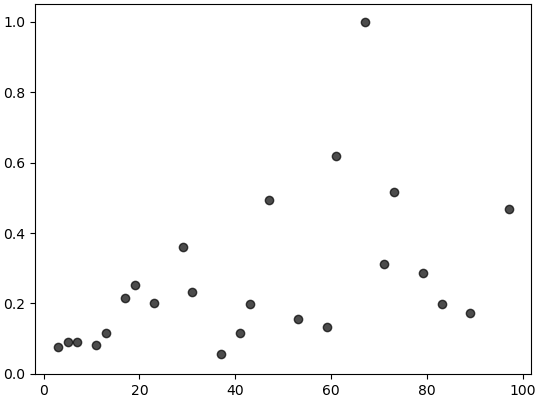}
        \caption{$p=2, \ell_p = 2, \ell_q=2$}
    \end{subfigure}

    \begin{subfigure}[b]{0.3\textwidth}
        \centering
        \includegraphics[width=\textwidth]{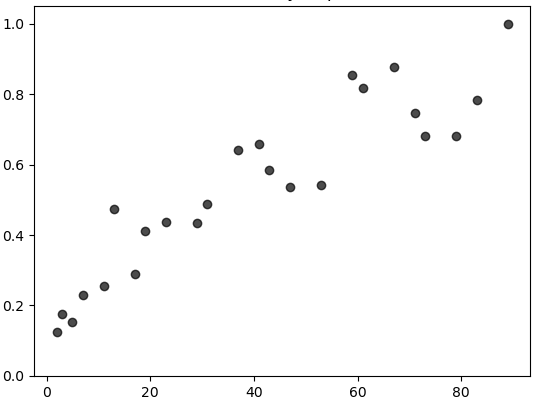}
        \caption{$p=97, \ell_p = 0, \ell_q=0$}
    \end{subfigure}
    \begin{subfigure}[b]{0.3\textwidth}
        \centering
        \includegraphics[width=\textwidth]{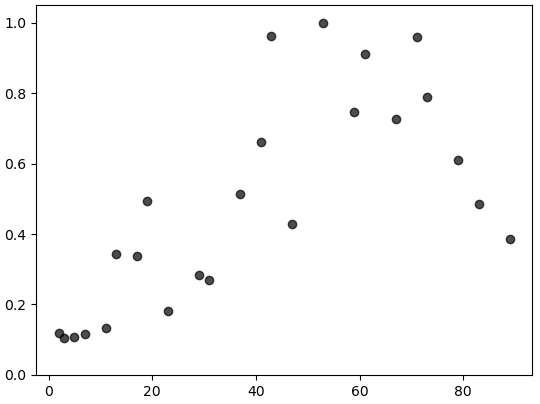}
        \caption{$p=97, \ell_p = 2, \ell_q=0$}
    \end{subfigure}
    \begin{subfigure}[b]{0.3\textwidth}
        \centering
        \includegraphics[width=\textwidth]{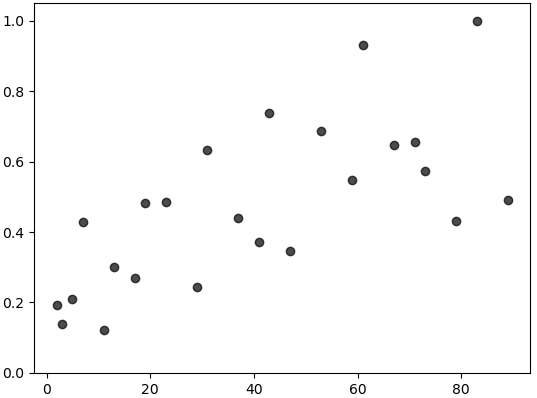}
        \caption{$p=97, \ell_p = 2, \ell_q=2$}
    \end{subfigure}

    \caption{Saliency plots (prime $q$ against saliency of $a_q \bmod \ell_q$) for neural network experiments in Table \ref{tab:NN_performance}}
    \label{fig:NN_saliency}
\end{figure}

One might expect that small primes carry more information about the $L$-function than larger primes. However, we obtain conflicting results, which require further investigation. We explore some observations on which $a_q$ appear to be most important for predicting $a_p$ by examining saliency plots. Figure \ref{fig:NN_saliency} shows saliency plots for the experiments in Table \ref{tab:NN_performance} in the same order as presented in the table.

We compare the saliency plot for the experiment for predicting $a_{97} \pmod{2}$ from subfigure (E) to five transformer-based experiments with architecture as described in \Cref{ss: int2int}. In each experiment, we exclude different sets of 3 primes. The models have the same encoder-only architecture and hyperparameters as the models in \Cref{p: ap}, except that we train the models for 150 epochs of size \numprint{100000} curves. We balance the datasets and use \numprint{1960662} curves, half of which have $a_{97} \equiv 0, 1 \pmod{2}$.

This comparison is appropriate since both the previous neural network experiment and the transformer initially treat the $a_q$ values as one-hot encodings. We record the maximum  MCC in Table \ref{tab: remove_primes}. We note that the experiments which exclude larger primes usually  have lower performance.

\begin{table}[!ht]
\centering
\begin{tabular}{|c|c|c|c|c|c|c|}
\hline
Primes excluded & Maximum Test MCC\\ \hline
2, 3, 5   &   0.8262 \\ \hline
13, 17, 19   & 0.7870 \\ \hline
31, 37, 41   &  0.7848 \\ \hline
43, 53, 71  &  0.7850\\ \hline
79, 83, 89  &  0.7658 \\ \hline
\end{tabular}
\caption{Maximum MCC values for predicting $a_{97}$ from $(a_q)_{q < 97}$ and excluding $a_q$ values for triples of primes. Compare with the saliency map in Subfigure (E) in Figure \ref{fig:NN_saliency}.
}
\label{tab: remove_primes}
\end{table}

\subsection{A different normalization}
We performed some experiments to predict $a_p \pmod 2$ using $\frac{a_q}{\sqrt{q}} < 2$ as input. This is in the same style as the experiments conducted in \ref{ss: a2NN}.

\begin{figure}[ht]
    \centering
    \begin{subfigure}{0.4\textwidth}
        \centering
        \includegraphics[width=\linewidth]{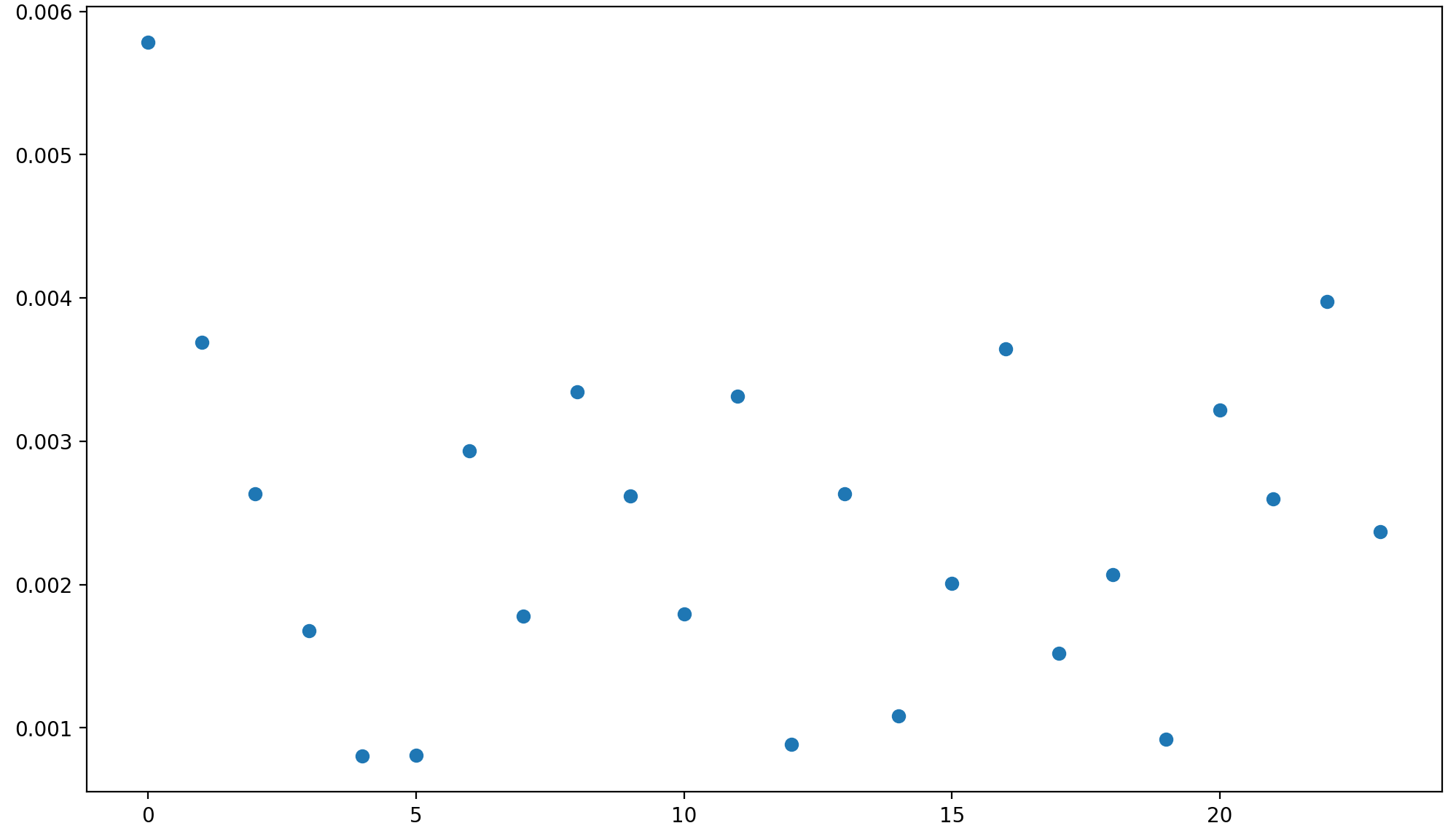}
    \end{subfigure}
    \hspace{0.1\textwidth}
    \begin{subfigure}{0.4\textwidth}
        \centering
        \includegraphics[width=\linewidth]{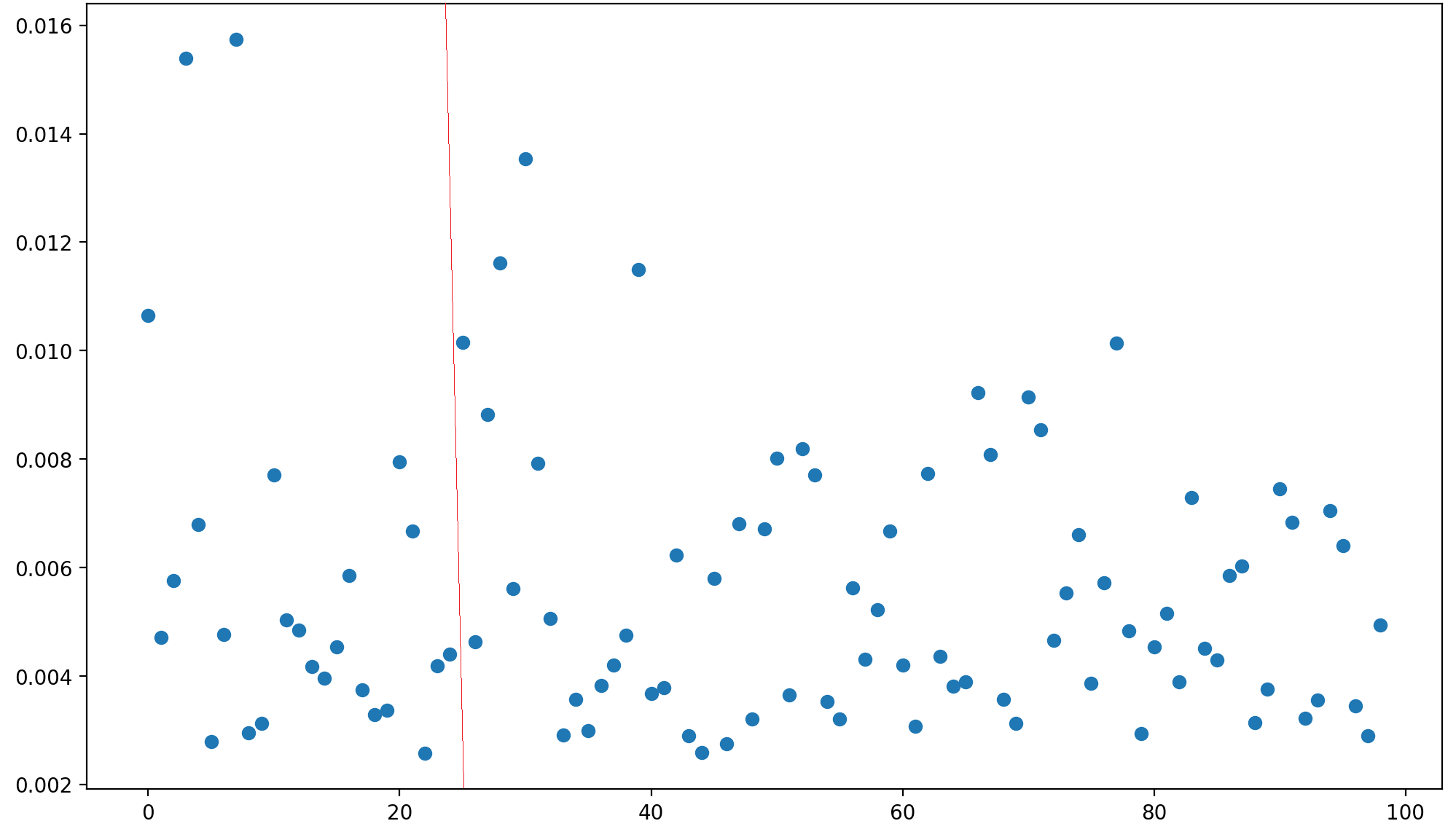}
    \end{subfigure}

    \vspace{0.5cm}
    \begin{subfigure}{0.4\textwidth}
        \centering
        \includegraphics[width=\linewidth]{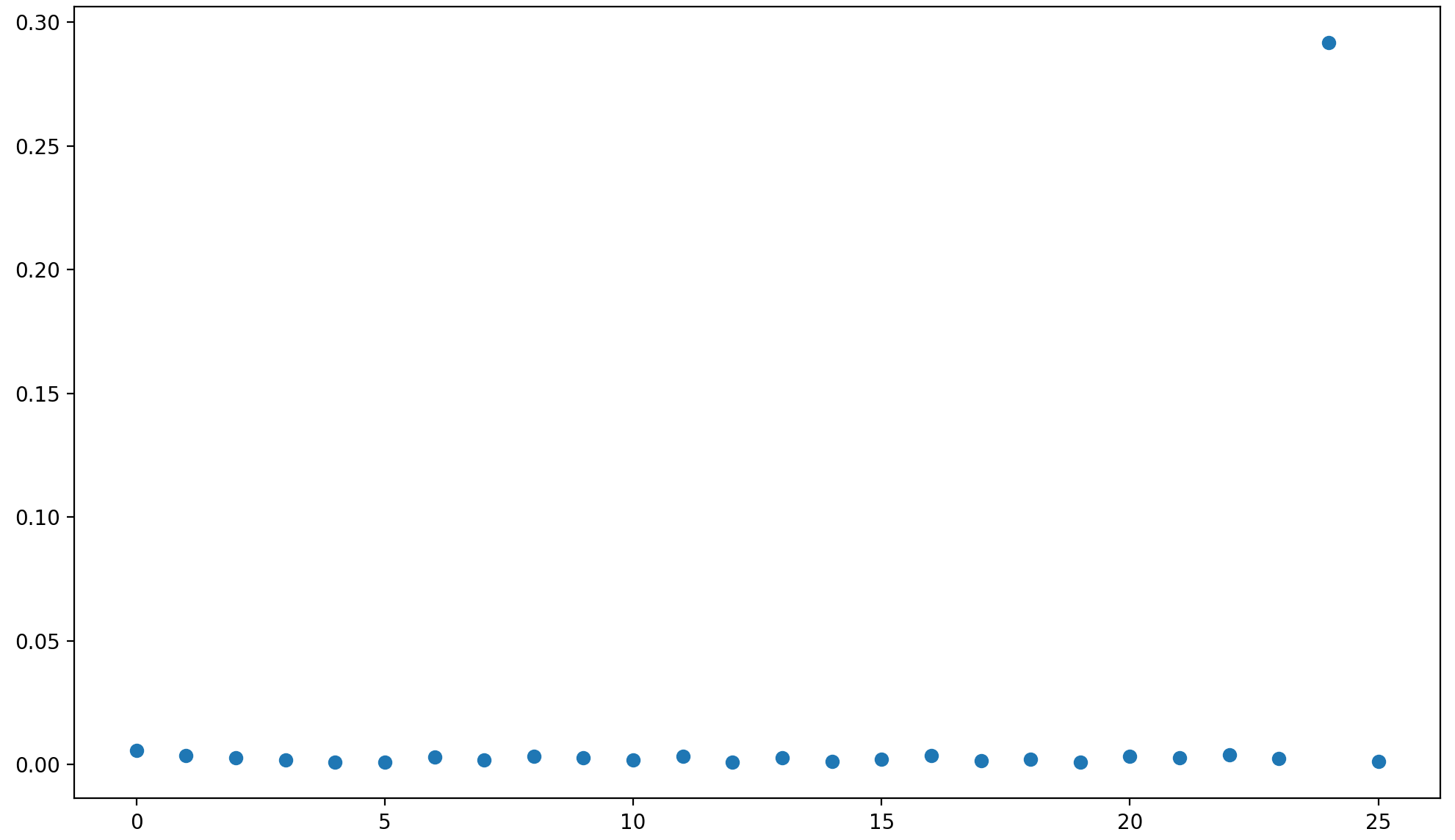}
    \end{subfigure}
    \hspace{0.1\textwidth}
    \begin{subfigure}{0.4\textwidth}
        \centering
        \includegraphics[width=\linewidth]{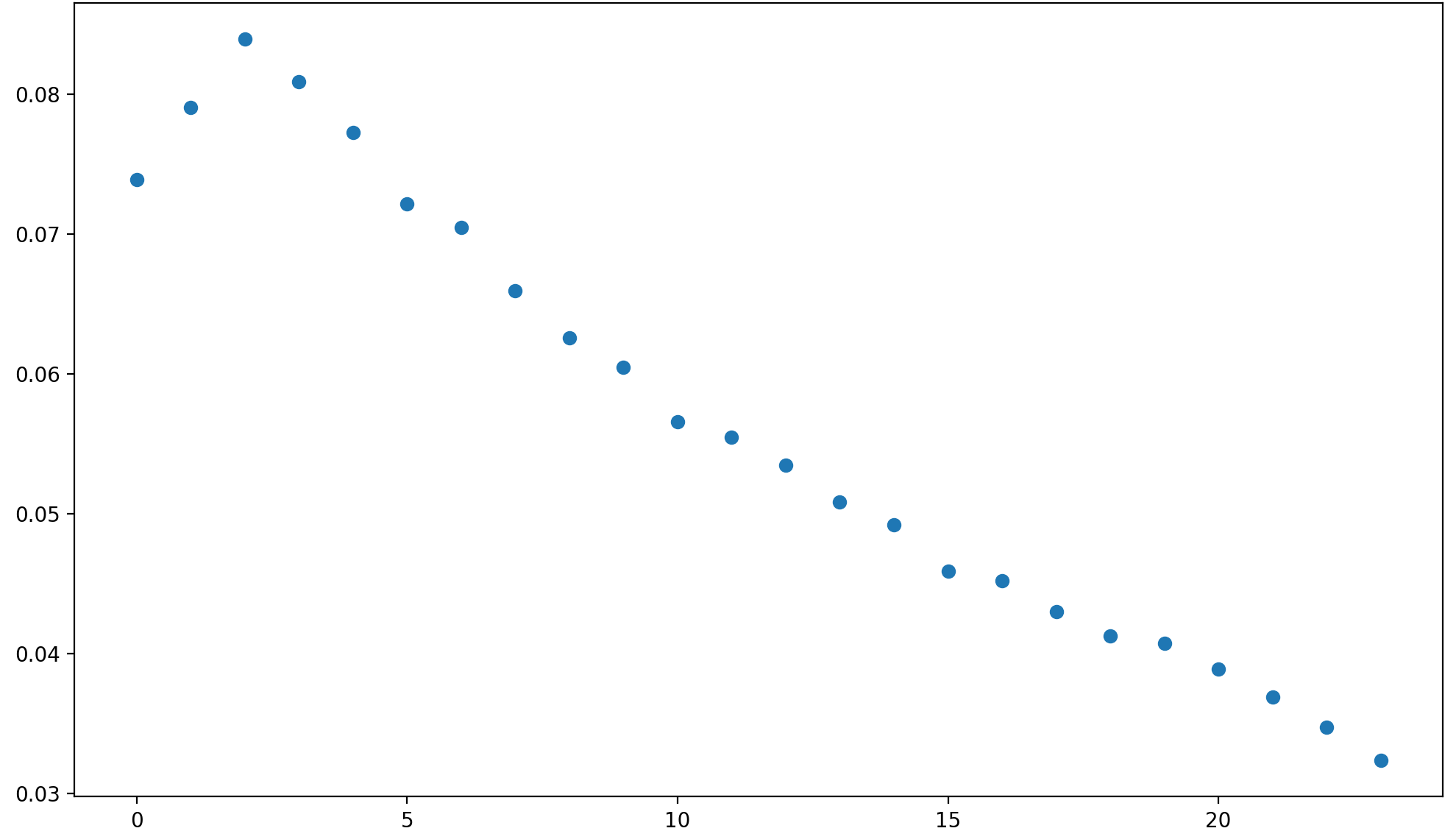}
    \end{subfigure}

    \caption{Saliency plots for important features. The top left is for the 24 primes governing $a_{97} \pmod 2$, the top right shows the importance of primes when using the first 100 primes. The bottom left shows the same plot with conductor and root number added as well, while the bottom right shows a saliency plot for $|a_2| = 1$.}
    \label{fig: a97_mod_2_ash}
\end{figure}

We use the dataset ECQ6, where we conducted experiments in the same style for a balanced dataset (size 1,000,000) with conductors between $10^5$ and $10^6$. For $p=97$, we obtained an accuracy of about 0.86, with an MCC of 0.7486, which is comparable to the performance of the FNNs in the previous section. The saliency plots (\ref{fig: a97_mod_2_ash}), in this case, however, do not show an independence from smaller primes. In fact, the small and middle primes seem to be the most important features. This is different to what we noticed with results from experiments using one-hot encoding, leading us to believe that a difference in style of encoding could lead to a difference in the importance of primes. Moreover, as expected, the conductor seems to be playing the most important role, while the root number has importance comparable to other prime features.

\subsection{Neural net experiments for \texorpdfstring{$a_2 \in \{ \pm 1\}$}{a2 in \{1, -1\}}}
\label{ss: a2NN} Initial models are trained on curves between conductors $2^{15}$ and $2^{16}$ (between $10^5$ and $10^6$) present in the LMFDB. For each curve in the dataset, $\frac{a_q}{\sqrt{q}} < 2$ values are given for the first 100 primes, except $a_2$, which is not normalized. The conductor and root number are also passed separately. \\
The total size of the balanced dataset is \numprint{71643}. Given the small size, we use 10\% of this as the test set. The training set is further scaled, while a test set is scaled using the same values. A one hot encoding is used to capture the output, and a final softmax layer ensures that the probabilities of the logits is returned. \\
A relatively small neural net with 2 hidden layers of dimensions 64 and 32 respectively gives the best result, with the AdamW optimizer having a dropout of 0.05, weight decay of 0.01 and learning rate $10^{-4}$. After \numprint{15000} epochs, we have an accuracy of 0.9137 on the test set.

We then repeated a similar experiment on the larger dataset ECQ6, which contains \numprint{1674207} datapoints with $a_2 \in \{\pm1\}$. Using the same parameters and replacing the final softmax layer with a sigmoid layer, after 250 epochs and with a batch size of \numprint{1024}, we obtain an accuracy of 0.6384 with an MCC of 0.3049. The striking thing about this model is its saliency graph, where the importance of the features seems to be linearly decreasing in nature, as shown in Figure \ref{fig: a97_mod_2_ash}.

\section{Visualizing Mod $\ell$ Structure in Learned Embeddings and Representations}
\label{p: interpret}

This section first examines the 2D PCA projections of the learned initial embeddings from the experiments conducted in \Cref{p: ap} and \Cref{p: apmod2}.
Our goal is to investigate whether the models in Section \ref{p: ap} predict \( a_p \) values by first learning to predict \( a_p \pmod{2} \) as an intermediate step.
Although tools for interpreting neural network-based models, including transformers, remain limited (for examples exploring how transformers learn some algorithmic tasks, see \cite{nanda2023progress, liu2022towards, zhong2024clock, quirke2023understanding, chughtai2023toy, hanna2024does, charton2023can}), the embeddings reveal clear patterns suggesting that the models encode distinctions based on their mod-\(\ell\) residues, with a particular emphasis on parity. The visualizations presented in this section were generated using  \cite{Bab}.

We further analyze the experimental results across several tasks: (1) predicting the exact values of \( a_p \) and (2) predicting their residues modulo \( \ell \) for \( \ell \in \{2, 4\} \). The model demonstrates significantly better performance in the \( a_p \pmod{2} \) and \( a_p \pmod{4} \) tasks.
Together with the PCA projections of the initial embeddings and the confusion matrix in \Cref{fig: a97_conf}, which reveals the model's difficulty in distinguishing values sharing the same parity, these results strongly suggest that the transformer model relies on predicting \( a_p \pmod{2} \) as an intermediate step in learning to predict the full \( a_p \) values.

Lastly, we implemented an additional encoder-decoder transformer model to predict \( a_2 \). The PCA projections of the decoder's hidden states and encoder's initial embeddings indicate that the model relies on inferences of \( a_p \pmod{\ell} \) for multiple \(\ell\).

\subsection{2D PCA Projections of Initial Embeddings}
The vocabularies used in the experiments described in Section \ref{p: ap} and \Cref{p: apmod2} consist of tokens including integers \( 0 \leq n \leq 19 \), the symbols \( + \) and \( - \), and special tokens like \texttt{<eos>}, which indicates the end of a sentence. Following \cite{elhage2021mathematical}, the initial step of learning involves the transformer model learning a context-free token embedding \( W_E \), which maps each token to a meaningful vector representation in \( \mathbb{R}^{256} \).

\subsubsection{Predicting \texorpdfstring{$a_p$}{ap}}
\begin{figure}[ht]
    \centering
    \begin{subfigure}{0.48\textwidth}
        \centering
        \includegraphics[width=\linewidth]{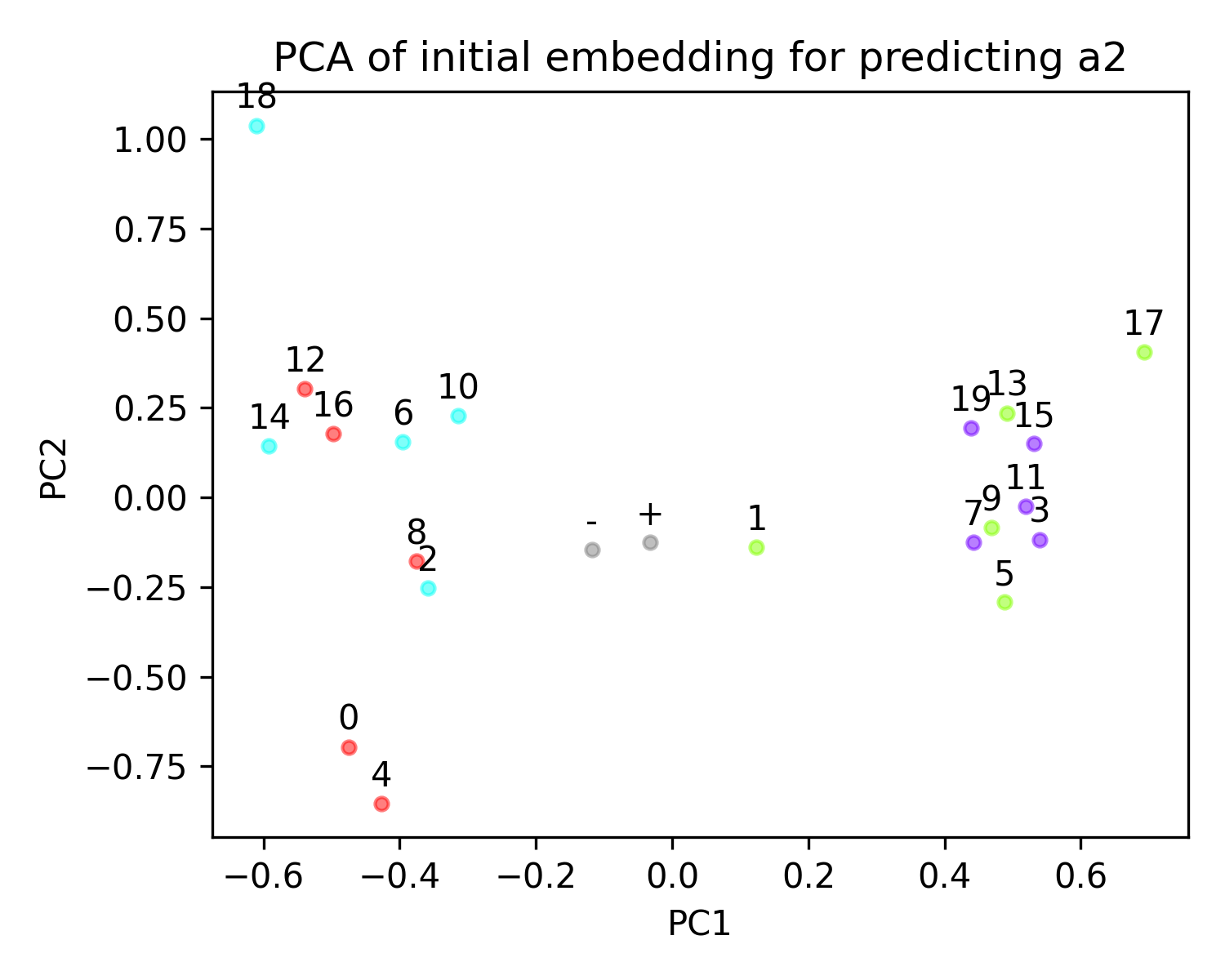}
    \end{subfigure}
    \begin{subfigure}{0.48\textwidth}
        \centering
        \includegraphics[width=\linewidth]{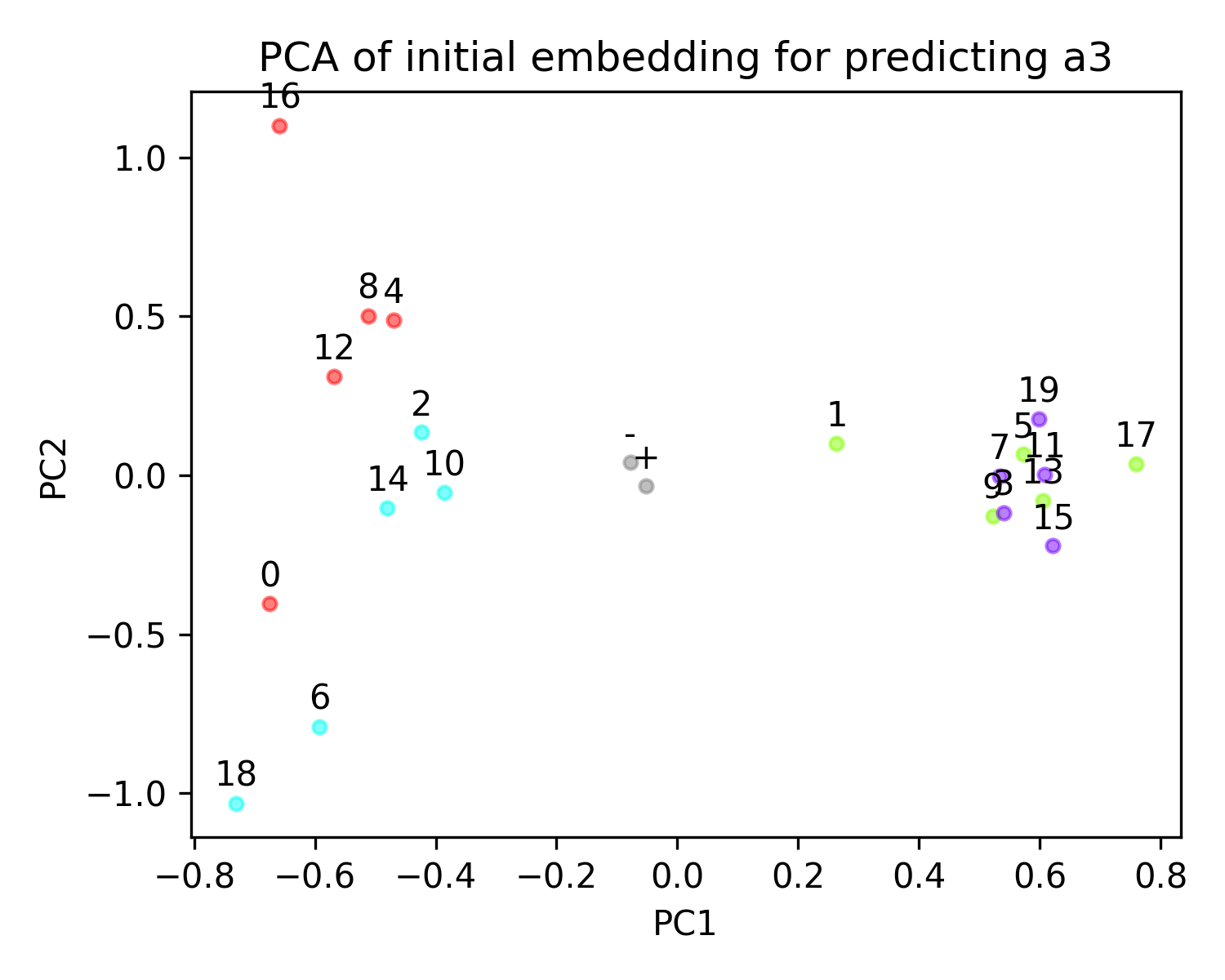}
    \end{subfigure}

    \vspace{0.5cm}
    \begin{subfigure}{0.48\textwidth}
        \centering
        \includegraphics[width=\linewidth]{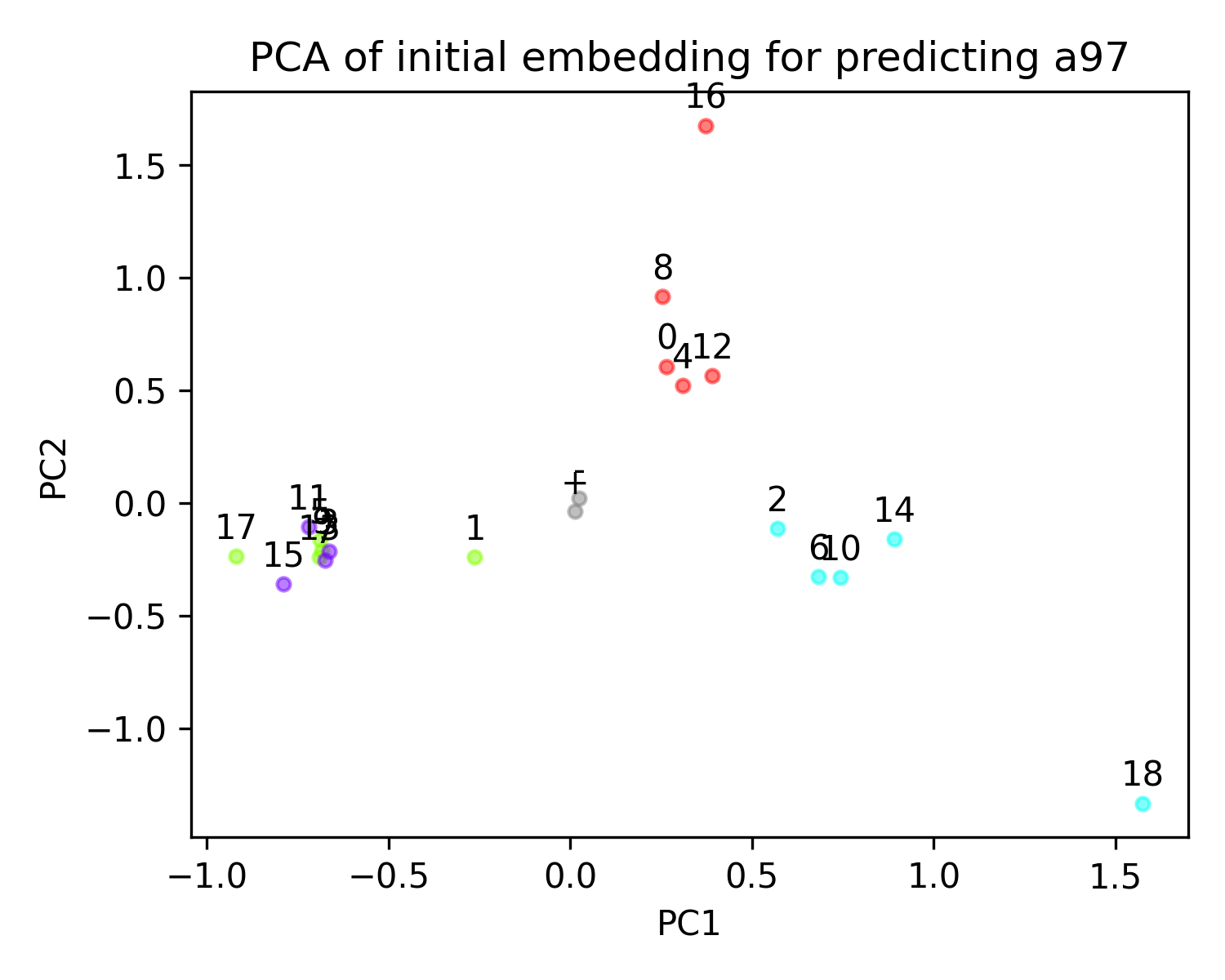}
    \end{subfigure}

    \caption{PCA plots of the initial embeddings for predicting \( a_2 \), \( a_3 \), and \( a_{97} \) are shown: \( a_2 \) (upper left), \( a_3 \) (upper right), and \( a_{97} \) (bottom).
    Red, green, blue, and purple dots represent \( a_q \equiv 0, 1, 2, 3 \pmod{4} \), respectively.
     All show clear separation between odd and even integers. For \( p \in \{3, 97\} \), we observe further differentiation between \( a_q \equiv 0 \pmod{4} \) and \( a_q \equiv 2 \pmod{4} \).
     }
    \label{fig: ap_initial_pca}
\end{figure}

{\label{sec: ap_pca}}
In \Cref{fig: ap_initial_pca}, we plot the 2D PCA projection of the initial embeddings for $+, -$ and integer values $0 \le n \le 19$ for the $a_p$ prediction task.

In \Cref{fig: ap_initial_pca}, there is a clear clustering based on the parity of the tokens in all cases, with further differentiation between \( a_q \equiv 0 \pmod{4} \) and \( a_q \equiv 2 \pmod{4} \) when predicting $a_3$ and $a_{97}$. These plots suggest that all the models implicitly learn to distinguish $a_q \pmod{2}$, and might therefore learn to predict $a_p \pmod{2}$ as an intermediate step.

In addition, in all cases, the \( + \) and \( - \) values cluster closely together near zero. As noted in Section \ref{p: ap}, this may suggest that the models have difficulty distinguishing between positive and negative \( a_p \) values and establishing a meaningful representation for them.

\subsubsection{Predicting \texorpdfstring{$a_p \pmod{2}$}{ap mod 2}} 
\begin{figure}[ht]
    \centering
    \begin{subfigure}{0.48\textwidth}
        \centering
        \includegraphics[width=\linewidth]{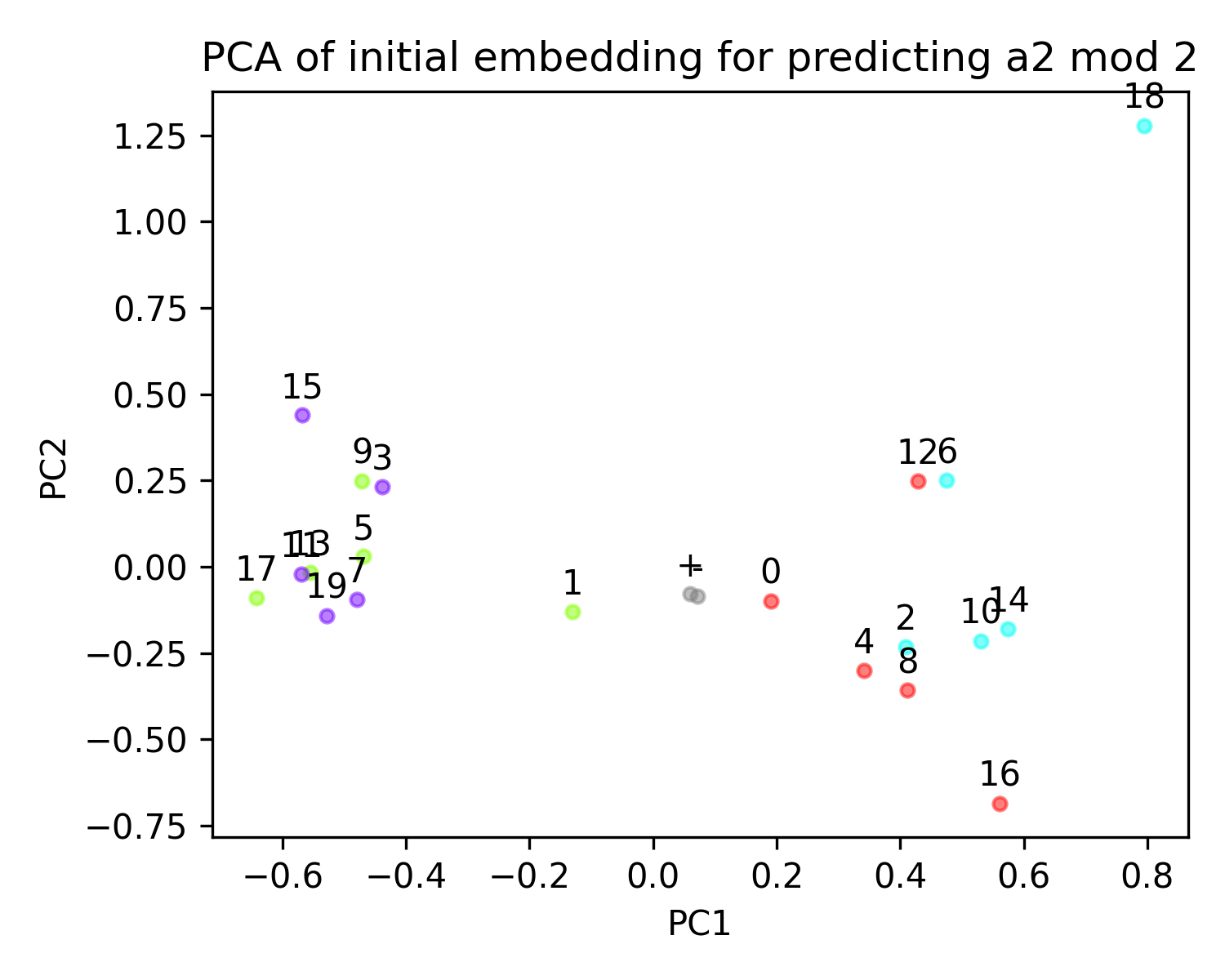}
    \end{subfigure}
    \begin{subfigure}{0.48\textwidth}
        \centering
        \includegraphics[width=\linewidth]{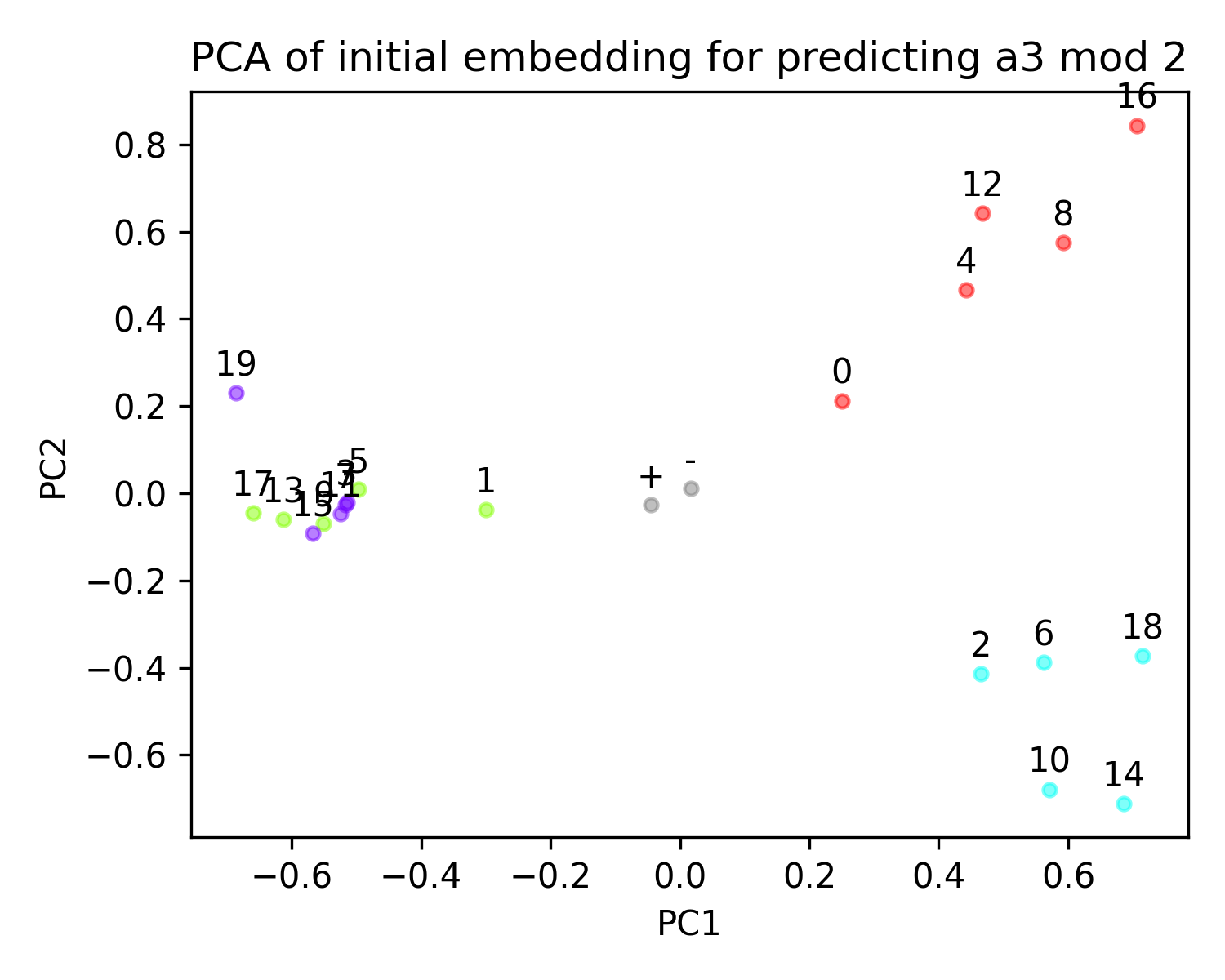}
    \end{subfigure}

    \vspace{0.5cm}
    \begin{subfigure}{0.48\textwidth}
        \centering
        \includegraphics[width=\linewidth]{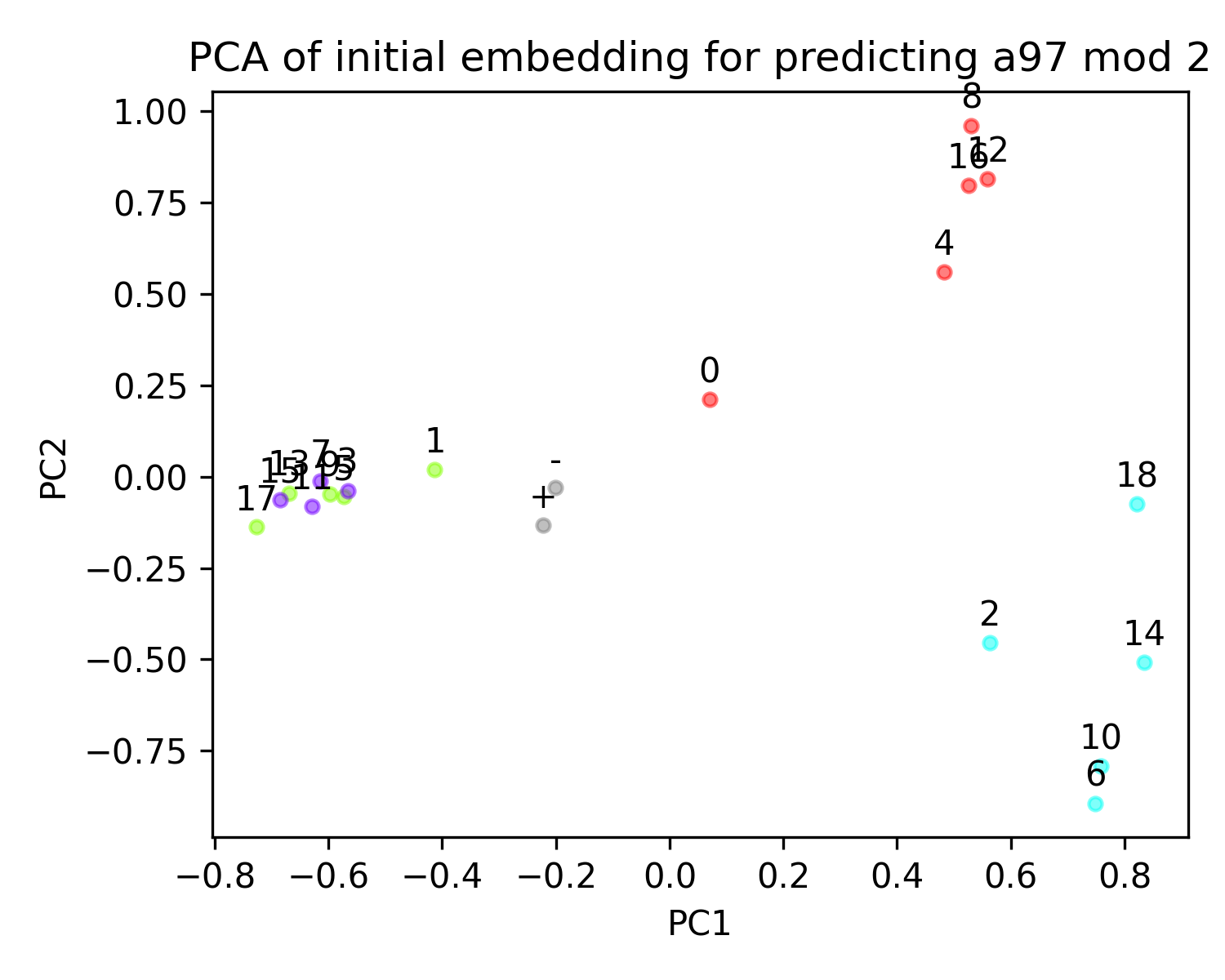}
    \end{subfigure}

    \caption{PCA plots of the initial embeddings for predicting  \( a_2 \pmod{2} \) (upper left), \( a_3 \pmod{2}\) (upper right), and \( a_{97} \pmod{2} \) (bottom).
    The distinction between modular--$4$ classes resembles that in \Cref{fig: ap_initial_pca}.
    }
    \label{fig: ap_mod2_pca}
\end{figure}

In \Cref{fig: ap_mod2_pca}, we plot the 2D PCA projection of the initial embeddings for the $a_p$ prediction task.  The differentiation between modular--$4$ classes is similar to \Cref{fig: ap_initial_pca} but slightly more pronounced.

\subsection{Experiment Results Comparison}

To further demonstrate that the transformer's learning replies on implicit predicting \( a_p \pmod{2} \), Table \ref{tab: learn_ap_mod} reports the MCCs for \( a_p \pmod{\ell} \) derived from the true \( a_p \) predictions in \Cref{p: ap}, for \( p \in \{2, 3, 97\} \) and \( \ell \in \{2, 3, 4\} \). These MCCs are computed by converting the labels to their modulo \( \ell \) values using the model trained on true \( a_p \) predictions. Notably, the performance for \( a_p \pmod{2} \) and \( a_p \pmod{4} \) predictions significantly surpasses that of predicting \( a_p \) and \( a_p \pmod{3} \), aligning with the model's recognition of mod-2 and mod-4 classes observed in \Cref{fig: ap_initial_pca}. Moreover, the \textit{Test MCC mod 2} results are consistent with, though slightly worse than, those in \Cref{p: apmod2}, where the model is specifically trained to predict \( a_p \pmod{2} \).

\begin{table}[!ht]
\centering
\begin{tabular}{|c|c|c|c|c|c|c|}
\hline
$p$ & $N_{\text{samples}}$ & $N_{\text{classes}}$ & Test MCC & Test MCC mod 2 & Test MCC mod 3 & Test MCC mod 4 \\ \hline
2   & \numprint{722204}        & 5        & 0.5822             & 0.7678 & 0.5385               & 0.6154     \\ \hline
3   & \numprint{1210417}        & 7        & 0.5205              & 0.7802    &   0.4893    & 0.5905    \\ \hline
97  & \numprint{3550727}        & 39        & 0.4711               & 0.7981      &0.4702        & 0.6497     \\ \hline
\end{tabular}
\caption{Predicting $a_p$ with an encoder-only model. The implied MCC results in \Cref{p: ap} are obtained by converting the labels to their modulo $\ell$ values and computing the MCC.
}
\label{tab: learn_ap_mod}
\end{table}

\subsection{An Encoder-Decoder Model}
For the purpose of visualizing decoder-specific learned embeddings and representations, we also trained an encoder-decoder transformer to predict \( a_2 \). The architecture and hyperparameters remain nearly identical to the previous models, with the only difference being the addition of a decoder containing one transformer layer. The architecture is depicted in \Cref{fig: encoder_decoder_archi}. Table \ref{tab: ap_dec_mod} reports the performance of this experiment.

\begin{table}[!ht]
\centering
\begin{tabular}{|c|c|c|c|c|c|c|c|}
\hline
$p$ & $N_{\text{samples}}$ & $N_{\text{classes}}$  & Test MCC & Test MCC mod 2 & Test MCC mod 3 & Test MCC mod 4 \\ \hline
2   & \numprint{722204}        & 5       & 0.6178            & 0.7788 & 0.5807              & 0.6485    \\ \hline
\end{tabular}
\caption{Predicting $a_2$ with an encoder-decoder model. The modular MCC results are obtained by converting the labels to their modulo $\ell$ values and computing the MCC.
}
\label{tab: ap_dec_mod}
\end{table}

\definecolor{DarkGray}{HTML}{68904d}
\definecolor{Orange}{HTML}{ffae42}
\begin{figure}[h!]
\small
\centering
\begin{tikzpicture}[
    box/.style={rectangle,draw,fill=DarkGray!20,node distance=1cm,text width=8em,text centered,rounded corners,minimum height=1.5em,thick},
    box2/.style={rectangle,draw,fill=DarkGray!20,node distance=1cm,text width=1em,text centered,rounded corners,minimum height=0.5em,thick},
    arrow/.style={draw,-latex',thick},
  ]
  \node [box] (ffne) {FCNN};

  \tikzset{arrow/.style={-latex, thick}}

  \node [box2,above=0.5 of ffne] (tok2) {c$_2$};
  \node [box2,left=0.25 of tok2] (tok1) {c$_1$};
  \node [box2,right=0.25 of tok2] (tok3) {c$_3$};
  \node [above=0.2 of tok2]{Encoder output};
  \node [box,below=0.23 of ffne,fill=Orange] (selfatt){Self-attention};  
  \node [label={[rotate=90]160:$\times n$},rectangle,draw,dashed,inner sep=0.7em,fit=(selfatt) (ffne)] (encoder1) {};

  \node [box2,below=0.4 of selfatt] (se){x$_2$};
  \node [box2,left=0.25 of se] (se2){x$_1$};
  \node [box2,right=0.25 of se] (se3){x$_3$};
  \node [box,below=0.23 of se] (embed) {Embedding};
  \node [below=0.3 of embed] (inp2) {\texttt{+ 2}};
\node [right=0.25 of inp2] (inp3) {\texttt{+ 5}};
\node [left=0.25 of inp2] (inp1) {\texttt{- 1}};
\path[arrow] (inp2) -- (embed);
\path[arrow] (inp3) -- (embed.south -| inp3);
\path[arrow] (inp1) -- (embed.south -| inp1);
\node [below=0.1 of inp2]{\textbf{Encoder}};

\node [box,right=1.0 of embed](emb2){Embedding};

\node [below=0.3 of emb2](enc2){\texttt{<BOS>}};
% \node [left=0.25 of enc2](enc3){\texttt{<BOS>}};
\node [below=0.1 of enc2]{\textbf{Decoder}};

\node [box2,above=0.23 of emb2](em22){y$_1$};
\path[draw,thick] (emb2) -- (em22);
\path[arrow] (enc2) -- (emb2);
% \path[arrow] (enc3) -- (emb2.south -| enc3);

% \node [box2,left=0.25 of em22](em21){y$_1$};
% \path[draw,thick] (emb2.north -| em21) -- (em21);

\node [box,above=0.4 of em22,fill=Orange] (selfatt2){Self-attention};
\path[draw,thick] (em22) -- (selfatt2);
% \path[draw,thick] (em21) -- (selfatt2.south -| em21);

% \node [box,above=0.23 of selfatt2] (gate){Copy-gate};
% \path[draw,thick] (selfatt2) -- (gate);

\node [box,above=0.23 of selfatt2,fill=Orange] (XattD){Cross-attention};
\path[draw,thick] (selfatt2) -- (XattD);

\node [box,above=0.23 of XattD] (ffn) {FCNN};
\path[draw,thick] (XattD) -- (ffn);

\node [label={[rotate=270]35:Layer 1},rectangle,draw,dashed,inner sep=0.9em,fit=(selfatt2) (ffn)] (decoder2) {};
% \node [box,above=0.55 of ffn,fill=Orange] (selfatt3){Self-attention};
% \path[draw,thick] (ffn) -- (selfatt3);

% \node [box,above=0.23 of selfatt3,fill=Orange] (crossatt2){Cross-attention};
% \path[draw,thick] (selfatt3) -- (crossatt2);

% \node [box,above=0.23 of crossatt2] (ffn2) {FCNN};
% \path[arrow] (crossatt2) -- (ffn2);
% \node [label={[rotate=270]22:Layer 2},rectangle,draw,dashed,inner sep=0.7em,fit=(selfatt3) (ffn2)] (decoder2) {};
\node [box,above=0.4 of ffn] (linear){Linear classifier};
\path[draw,thick] (ffn) -- (linear);
\node [box,above=0.23 of linear] (softmax){Softmax};
\path[draw,thick] (linear) -- (softmax);
\node[above=0.3 of softmax](output){\texttt{-1}};
\path[arrow] (softmax) -- (output);

\coordinate[above right =0.2 and 0.7 of tok3](conn);
\path[draw,thick] (tok1) |- (conn);
\path[draw,thick] (tok2) |- (conn);
\path[draw,thick] (tok3) |- (conn);

% \coordinate[left=0.1 of gate](lg);
% \path[draw,thick] (gate) -- (lg);
\coordinate[above=0.1 of ffn.north](mw2);
% \path[draw,thick] (lg) |- (mw2);
\coordinate[right=1.7 of mw2](mw3);
\path[draw,thick] (mw2)-- (mw3);
% \coordinate[below right =0.1 and 0.1 of selfatt2.south](mw);
\coordinate[below =0.1 of selfatt2.south](mw);X
\path[draw,thick] (mw3) |- (mw);
 
\path [draw,thick] (ffne.north -| tok1) -- (tok1);
\path [draw,thick] (ffne) -- (tok2);
\path [draw,thick] (ffne.north -| tok3) -- (tok3);
\path [draw,thick] (selfatt) -- (ffne);
\path [draw,thick] (se) -- (selfatt);
\path [draw,thick] (se2) -- (selfatt.south -| se2);
\path [draw,thick] (se3) -- (selfatt.south -| se3);
\path[draw,thick] (embed.north -| se2) -- (se2);
\path[draw,thick] (embed) -- (se);
\path[draw,thick] (embed.north -| se3) -- (se3);
% \coordinate[left=0.1 of XattD.west](mw4);
\path[draw,thick] (conn) |- (XattD);
% \path[draw,thick] (conn) |- (crossatt2.west);

\node[draw,text width=5cm] at (-5.2,0) {Embedding dimension\\
\begin{itemize}[leftmargin=*, itemsep=1pt, topsep=1pt, parsep=1pt]
\item Encoder: $d_e=256$ \\
\item Decoder: $d_e=256$ 
\end{itemize}
Attention heads: $8$

Fully connected neural network (FCNN):\\ 
\begin{itemize}[leftmargin=*, itemsep=1pt, topsep=1pt, parsep=1pt]
\item 1 hidden layer \\
\item $4d_e$ neurons \\
\end{itemize}
Encoder Transformer layer:\\
  \begin{itemize}[leftmargin=*, itemsep=1pt, topsep=1pt, parsep=1pt]
      \item Iterated: $n = 4$ times
  \end{itemize}};
\end{tikzpicture}
\vspace{-0.1cm}
\caption{We also implement an encoder-decoder Transformer for further insights into $a_2$ predictions using Int2Int.} \label{fig: encoder_decoder_archi}
%\vspace{-0.1cm}
\end{figure}
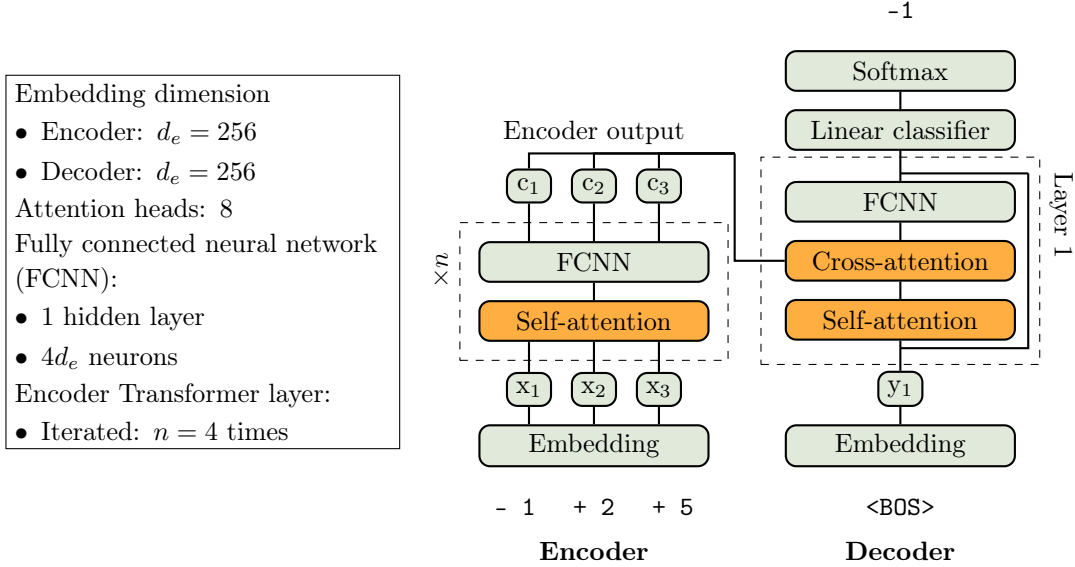

\begin{figure}[htb!]
    \centering
        \includegraphics[width=0.55\textwidth, height=7cm]{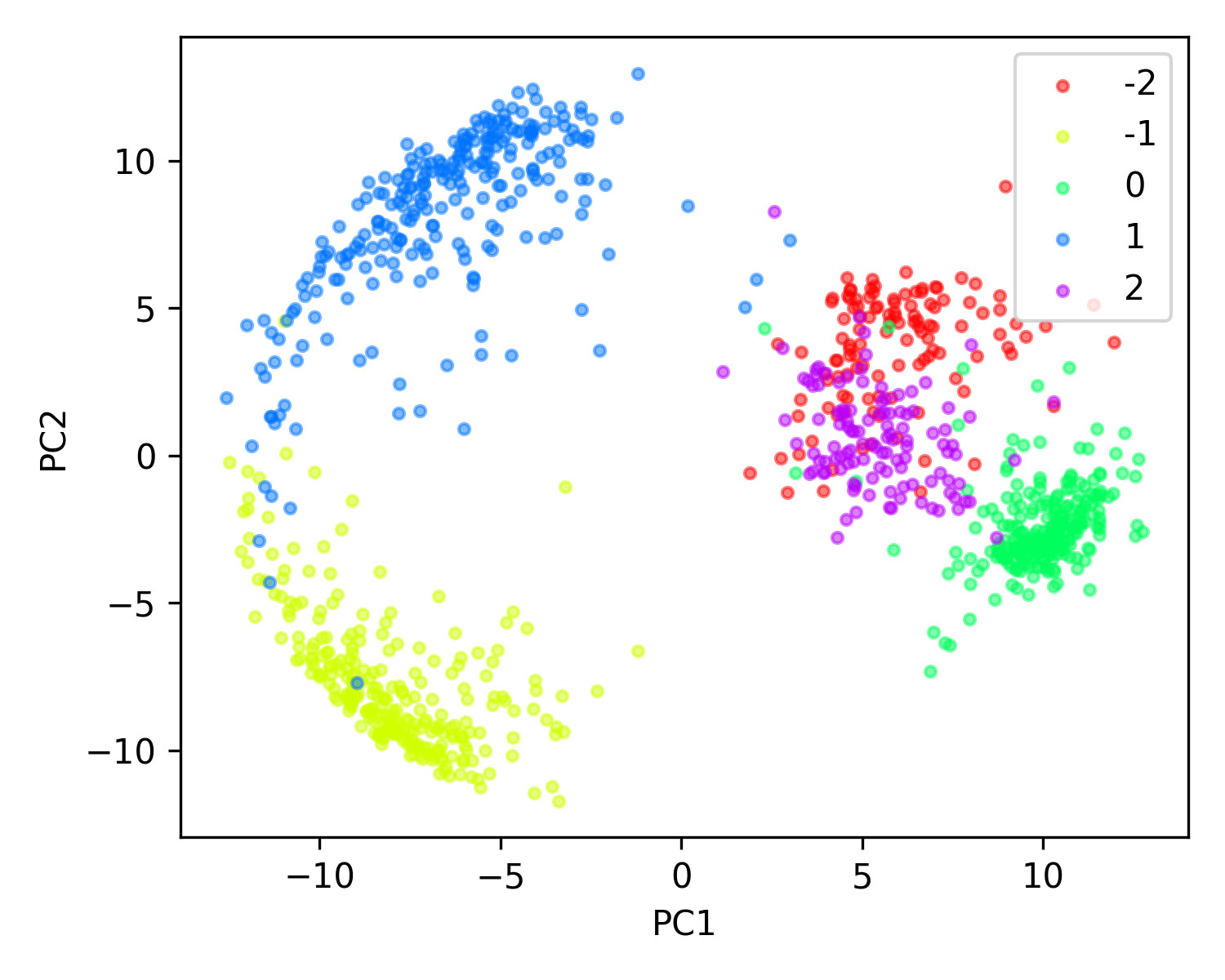}
        \caption{The 2D PCA of the model’s hidden states at the decoder transformer layer for predicting $a_2$ with conductor up to $10^4$. The dots, representing the hidden states of curves, are separated by their residue classes mod-2 and mod-4 from left to right.
        }
         \label{fig: a2_decodePCA}
\end{figure}

\begin{figure}[htb!]
    \centering
    \begin{subfigure}{0.48\textwidth}
        \centering
        \includegraphics[width=\linewidth]{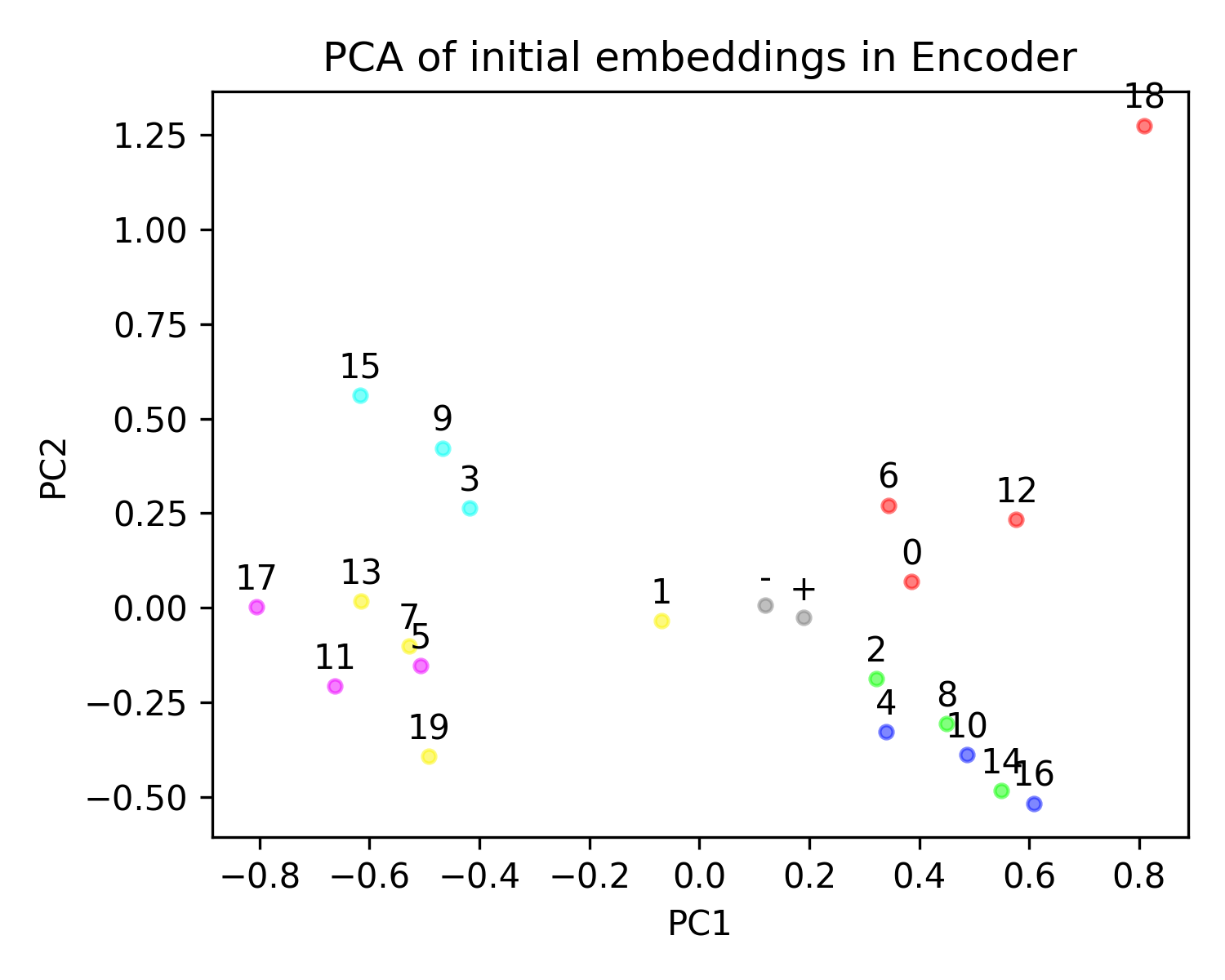}
    \end{subfigure}
    \begin{subfigure}{0.48\textwidth}
        \centering
        \includegraphics[width=\linewidth]{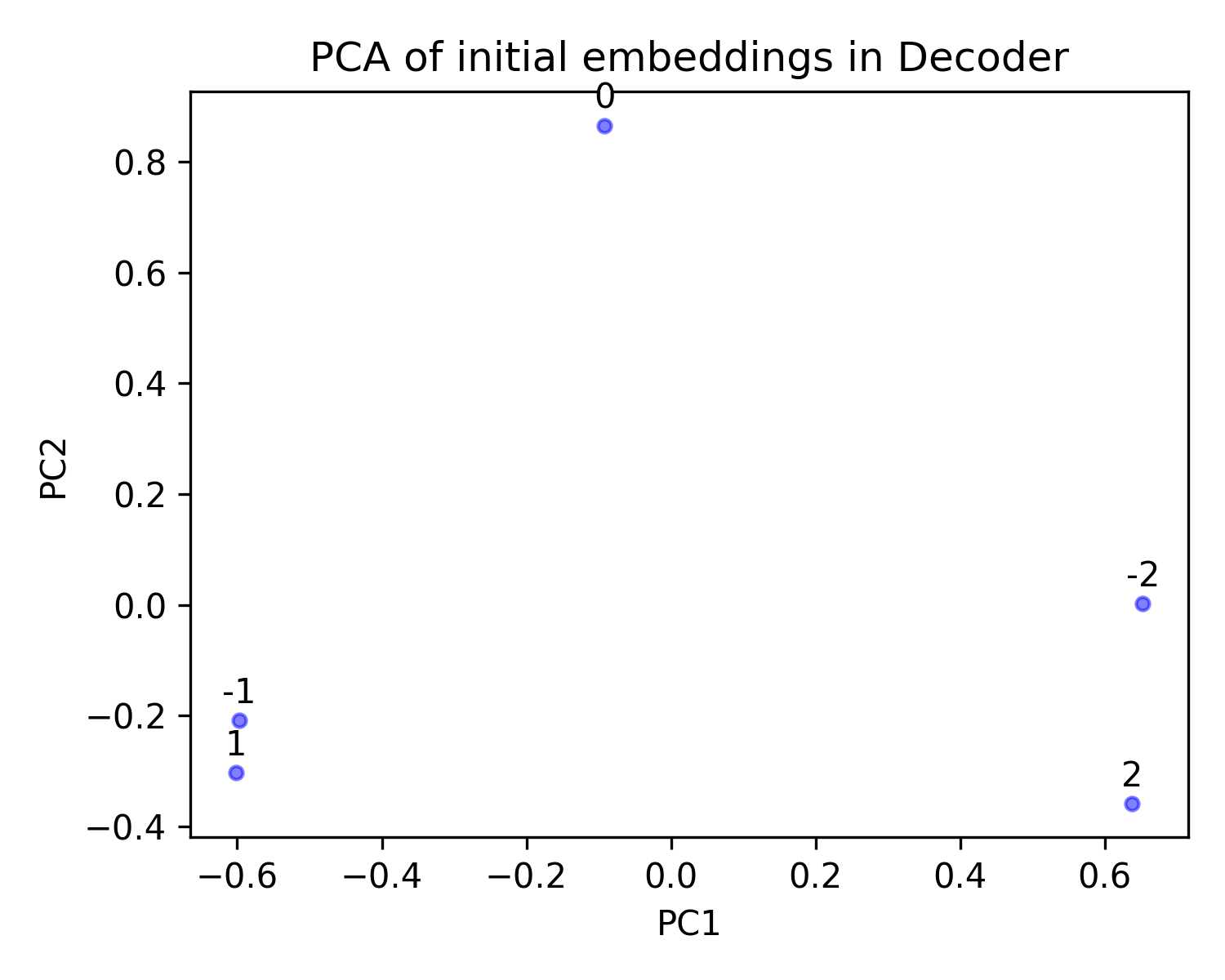}
    \end{subfigure}

    \caption{PCA plots of the initial embeddings in the encoder (left) and decoder (right) for predicting \( a_2 \) with an encoder-decoder model. The encoder PCA shows a clear separation between even and odd integers, as well as integers divisible by 3. As a result, we observe clustering based on modulo 6.
     }
    \label{fig: ap_initial_pca_enc_dec}
\end{figure}

As a result, we obtain a clear structure of the hidden states in the decoder compared to the encoder. Figure \ref{fig: a2_decodePCA} shows the 2D PCA plot of the decoder's hidden states after passing the model curves of low conductor up to $10^4$. The model distinguishes these elliptic curves by their \( a_2 \pmod{2} \) and \( a_2 \pmod{4} \) values after processing the \( a_q \) inputs. Additionally, the second principal component reflects a distinction between positive and negative values, possibly due to the model's recognition of \( a_2 \pmod{3} \) values.

We also note differences in the learned initial embedding in the encoder, as shown in Figure \ref{fig: ap_initial_pca_enc_dec}. Compared to the encoder-only model used in \Cref{p: ap} visualized in Figure \ref{fig: ap_initial_pca}, the encoder's initial embeddings in this experiment contain clusters based both on parity and divisibility by 3, resulting in clustering based modulo 6. Despite this clustering, the prediction of $a_2 \pmod{3}$ does not outperform the prediction of $a_2$.

\bibliographystyle{alpha}
\bibliography{bibfile}
\end{document}